\documentclass[12pt]{article}
\usepackage{amsfonts}
\usepackage{full page,
amssymb, amscd, amsmath,graphicx, 
}
\usepackage{amsthm}
\usepackage[margin=1in]{geometry}
 \title{{\bf Vertex algebraic intertwining operators among generalized Verma
modules for $\widehat{\mathfrak{sl}(2,\mathbb{C})}$}}
 \author{Robert McRae and Jinwei Yang}
    \date{}
    \begin{document}
    \bibliographystyle{alpha}
    \maketitle
\newtheorem{thm}{Theorem}[section]
\newtheorem{defn}[thm]{Definition}
\newtheorem{prop}[thm]{Proposition}
\newtheorem{cor}[thm]{Corollary}
\newtheorem{lemma}[thm]{Lemma}
\theoremstyle{definition}\newtheorem{rema}[thm]{Remark}
\newtheorem{app}[thm]{Application}
\newtheorem{prob}[thm]{Problem}
\newtheorem{conv}[thm]{Convention}
\newtheorem{conj}[thm]{Conjecture}
\newtheorem{cond}[thm]{Condition}
    \newtheorem{exam}[thm]{Example}
\newtheorem{assum}[thm]{Assumption}
     \newtheorem{nota}[thm]{Notation}
\newcommand{\halmos}{\rule{1ex}{1.4ex}}
\newcommand{\pfbox}{\hspace*{\fill}\mbox{$\halmos$}}

\newcommand{\nn}{\nonumber \\}

 \newcommand{\res}{\mbox{\rm Res}}
 \newcommand{\e}{\mbox{\rm exp}}
 \newcommand{\ord}{\mbox{\rm ord}}
\renewcommand{\hom}{\mbox{\rm Hom}}
\newcommand{\edo}{\mbox{\rm End}\ }
 \newcommand{\pf}{{\it Proof.}\hspace{2ex}}
 \newcommand{\epf}{\hspace*{\fill}\mbox{$\halmos$}}
 \newcommand{\epfv}{\hspace*{\fill}\mbox{$\halmos$}\vspace{1em}}
 \newcommand{\epfe}{\hspace{2em}\halmos}
\newcommand{\nord}{\mbox{\scriptsize ${\circ\atop\circ}$}}
\newcommand{\wt}{\mbox{\rm wt}\ }
\newcommand{\swt}{\mbox{\rm {\scriptsize wt}}\ }
\newcommand{\lwt}{\mbox{\rm wt}^{L}\;}
\newcommand{\rwt}{\mbox{\rm wt}^{R}\;}
\newcommand{\slwt}{\mbox{\rm {\scriptsize wt}}^{L}\,}
\newcommand{\srwt}{\mbox{\rm {\scriptsize wt}}^{R}\,}
\newcommand{\clr}{\mbox{\rm clr}\ }
\newcommand{\tr}{\mbox{\rm Tr}}
\newcommand{\C}{\mathbb{C}}
\newcommand{\Z}{\mathbb{Z}}
\newcommand{\R}{\mathbb{R}}
\newcommand{\Q}{\mathbb{Q}}
\newcommand{\N}{\mathbb{N}}
\newcommand{\CN}{\mathcal{N}}
\newcommand{\F}{\mathcal{F}}
\newcommand{\I}{\mathcal{I}}
\newcommand{\V}{\mathcal{V}}
\newcommand{\one}{\mathbf{1}}
\newcommand{\BY}{\mathbb{Y}}
\newcommand{\ds}{\displaystyle}
\newcommand{\gvmzero}{V_{\widehat{\mathfrak{g}}}(\ell,0)}
\newcommand{\g}{\mathfrak{g}}
\newcommand{\ghat}{\widehat{\mathfrak{g}}}
\newcommand{\gtilde}{\widetilde{\mathfrak{g}}}

        \newcommand{\ba}{\begin{array}}
        \newcommand{\ea}{\end{array}}
        \newcommand{\be}{\begin{equation}}
        \newcommand{\ee}{\end{equation}}
        \newcommand{\bea}{\begin{eqnarray}}
        \newcommand{\eea}{\end{eqnarray}}
         \newcommand{\lbar}{\bigg\vert}
        \newcommand{\p}{\partial}
        \newcommand{\dps}{\displaystyle}
        \newcommand{\bra}{\langle}
        \newcommand{\ket}{\rangle}

        \newcommand{\ob}{{\rm ob}\,}
        \renewcommand{\hom}{{\rm Hom}}

\newcommand{\A}{\mathcal{A}}
\newcommand{\Y}{\mathcal{Y}}

\newcommand{\dlt}[3]{#1 ^{-1}\delta \bigg( \frac{#2 #3 }{#1 }\bigg) }

\newcommand{\dlti}[3]{#1 \delta \bigg( \frac{#2 #3 }{#1 ^{-1}}\bigg) }

 \makeatletter
\newlength{\@pxlwd} \newlength{\@rulewd} \newlength{\@pxlht}
\catcode`.=\active \catcode`B=\active \catcode`:=\active
\catcode`|=\active
\def\sprite#1(#2,#3)[#4,#5]{
   \edef\@sprbox{\expandafter\@cdr\string#1\@nil @box}
   \expandafter\newsavebox\csname\@sprbox\endcsname
   \edef#1{\expandafter\usebox\csname\@sprbox\endcsname}
   \expandafter\setbox\csname\@sprbox\endcsname =\hbox\bgroup
   \vbox\bgroup
  \catcode`.=\active\catcode`B=\active\catcode`:=\active\catcode`|=\active
      \@pxlwd=#4 \divide\@pxlwd by #3 \@rulewd=\@pxlwd
      \@pxlht=#5 \divide\@pxlht by #2
      \def .{\hskip \@pxlwd \ignorespaces}
      \def B{\@ifnextchar B{\advance\@rulewd by \@pxlwd}{\vrule
         height \@pxlht width \@rulewd depth 0 pt \@rulewd=\@pxlwd}}
      \def :{\hbox\bgroup\vrule height \@pxlht width 0pt depth
0pt\ignorespaces}
      \def |{\vrule height \@pxlht width 0pt depth 0pt\egroup
         \prevdepth= -1000 pt}
   }
\def\endsprite{\egroup\egroup}
\catcode`.=12 \catcode`B=11 \catcode`:=12 \catcode`|=12\relax
\makeatother

\def\hboxtr{\FormOfHboxtr} 
\sprite{\FormOfHboxtr}(25,25)[0.5 em, 1.2 ex] 

:BBBBBBBBBBBBBBBBBBBBBBBBB | :BB......................B |
:B.B.....................B | :B..B....................B |
:B...B...................B | :B....B..................B |
:B.....B.................B | :B......B................B |
:B.......B...............B | :B........B..............B |
:B.........B.............B | :B..........B............B |
:B...........B...........B | :B............B..........B |
:B.............B.........B | :B..............B........B |
:B...............B.......B | :B................B......B |
:B.................B.....B | :B..................B....B |
:B...................B...B | :B....................B..B |
:B.....................B.B | :B......................BB |
:BBBBBBBBBBBBBBBBBBBBBBBBB |

\endsprite
\def\shboxtr{\FormOfShboxtr} 
\sprite{\FormOfShboxtr}(25,25)[0.3 em, 0.72 ex] 

:BBBBBBBBBBBBBBBBBBBBBBBBB | :BB......................B |
:B.B.....................B | :B..B....................B |
:B...B...................B | :B....B..................B |
:B.....B.................B | :B......B................B |
:B.......B...............B | :B........B..............B |
:B.........B.............B | :B..........B............B |
:B...........B...........B | :B............B..........B |
:B.............B.........B | :B..............B........B |
:B...............B.......B | :B................B......B |
:B.................B.....B | :B..................B....B |
:B...................B...B | :B....................B..B |
:B.....................B.B | :B......................BB |
:BBBBBBBBBBBBBBBBBBBBBBBBB |

\endsprite

\vspace{2em}

\begin{abstract}
\noindent We construct vertex algebraic intertwining operators among certain
generalized Verma modules for $\widehat{\mathfrak{sl}(2,\C)}$ and calculate the
corresponding fusion rules. Additionally, we show that under some conditions
these intertwining operators descend to intertwining operators among one
generalized Verma module and two (generally non-standard) irreducible modules.
Our construction relies on the irreducibility of the maximal proper submodules
of generalized Verma modules appearing in the Garland-Lepowsky resolutions of
standard $\widehat{\mathfrak{sl}(2,\C)}$-modules. We prove this irreducibility
using the composition factor multiplicities of irreducible modules in Verma modules for symmetrizable Kac-Moody Lie algebras of rank $2$, given by Rocha-Caridi and Wallach.
\end{abstract}



\renewcommand{\theequation}{\thesection.\arabic{equation}}
\renewcommand{\thethm}{\thesection.\arabic{thm}}
\setcounter{equation}{0} \setcounter{thm}{0}
\date{}
\maketitle

\section{Introduction}

Intertwining operators are fundamental objects in the representation theory of vertex operator algebras and in the construction of conformal field theories. Intertwining operators among a triple of modules for a vertex operator
algebra $V$ correspond to $V$-module homomorphisms from the tensor product of
two modules into the third (\cite{HLZ1}, \cite{HL}). Since tensor products of
$V$-modules cannot be obtained from the tensor products of the underlying vector
spaces, intertwining operators are essential to understanding the representation
theory of $V$. Indeed, the existence of tensor products depends on the
finiteness of the fusion rules, that is, the dimensions of spaces of intertwining
operators. Moreover, associativity of intertwining operators (\cite{H1},
\cite{HLZ2}) and modular invariance for traces of products of intertwining
operators (\cite{H2}) are crucial for showing that modules for certain vertex
operator algebras form modular tensor categories (\cite{H4}; see also the
review article \cite{HL}).

In \cite{FZ}, for vertex operator algebras for which a suitable category of weak modules is semisimple, Frenkel and Zhu identified the spaces of
intertwining operators among $V$-modules with suitable spaces constructed from
left modules and bimodules for the Zhu's algebra $A(V)$. In \cite{Li3}, Li gave a generalization and a proof of this result. In particular, this approach can be used to describe the intertwining operators among standard (that is, integrable highest weight) modules for the affine Lie algebra $\widehat{\mathfrak{g}}$ (where $\mathfrak{g}$ is a finite-dimensional simple Lie algebra) at a fixed level $\ell\in\N$. These $\widehat{\mathfrak{g}}$-modules are modules for the generalized Verma module vertex operator algebra $V^{M(\ell\Lambda_0)}$, where $\Lambda_0$ is the fundamental weight of $\widehat{\mathfrak{g}}$ corresponding the affine simple root $\alpha_0$.

In the present paper, we construct intertwining operators among generalized Verma modules for $V^{M(\ell\Lambda_0)}$ and their (generally non-standard) irreducible quotients in the case that $\mathfrak{g} = \mathfrak{sl}(2, \C)$. Since generalized Verma modules are typically reducible but indecomposable, the category of
$V^{M(\ell\Lambda_0)}$-modules is not semisimple, and we cannot use the method in \cite{FZ} for constructing intertwining operators among generalized Verma modules for $V^{M(\ell\Lambda_0)}$ or for calculating fusion rules. Instead, we develop new methods for constructing intertwining operators, which use the structure of the generalized Verma modules under consideration.

We use two results from the representation theory of Kac-Moody Lie algebras to
determine the structure of generalized Verma modules for $\widehat{\mathfrak{sl}(2,\C)}$. The first result, which more generally illustrates the
importance of generalized Verma modules in representation theory, is the
construction by Garland and Lepowsky in \cite{GL} of resolutions of standard
modules for symmetrizable Kac-Moody Lie algebras. These resolutions enabled
Garland and Lepowsky to prove the Macdonald-Kac formulas for symmetrizable
Kac-Moody Lie algebras. For $\widehat{\mathfrak{sl}(2,\C)}$, we focus on resolutions by generalized Verma modules: If $\Lambda$ is a dominant
integral weight of $\widehat{\mathfrak{sl}(2,\C)}$, there is an exact sequence
\[\cdots \rightarrow V^{M(w_j(\Lambda + \rho) - \rho)}
\stackrel{d_j}{\rightarrow}
V^{M(w_{j-1}(\Lambda + \rho) -\rho)}
\stackrel{d_{j-1}}{\rightarrow} \cdots
\stackrel{d_2}{\rightarrow} V^{M(w_1(\Lambda +\rho) - \rho)}
\stackrel{d_1}{\rightarrow} V^{M(\Lambda)}
\stackrel{\Pi_{\Lambda}}{\rightarrow} L(\Lambda) \rightarrow
0,\]
where the $d_j$ are $\widehat{\mathfrak{sl}(2,\C)}$-module
homomorphisms, $\Pi_\Lambda$ denotes the projection from $V^{M(\Lambda)}$ to its
irreducible quotient $L(\Lambda)$, and the $w_j$ are certain elements of length $j$ in
the Weyl group.

Secondly, we use a multiplicity result of Rocha-Caridi and Wallach \cite{RW2} for symmetrizable Kac-Moody Lie algebras of rank $2$, that composition factor multiplicities of irreducible modules in Verma modules equal $1$ or $0$. This is a special case of the Kazhdan-Lusztig multiplicity formula, conjectured in \cite{DGK} and proved by Kashiwara-Tanisaki in \cite{Ka}, \cite{KaT} and Casian in \cite{C}, that composition factor multiplicities are given by values of Kazhdan-Lusztig polynomials. Combining the result of Rocha-Caridi and Wallach with the Garland-Lepowsky resolutions, we prove our main result on the structure of generalized Verma modules for $\widehat{\mathfrak{sl}(2,\C)}$:  for a generalized Verma
$\widehat{\mathfrak{sl}(2,\C)}$-module appearing in the Garland-Lepowsky resolution of
a standard $\widehat{\mathfrak{sl}(2,\C)}$-module, the maximal proper submodule is irreducible.

Our main result states that under certain conditions, there is a natural vector space isomorphism between the space of
intertwining operators of type $\binom{V^{M(r)}}{V^{M(p)}\,V^{M(q)}}$ and
\begin{equation*}
 \mathrm{Hom}_{\mathfrak{sl}(2,\C)}(M(p)\otimes M(q),M(r)),
\end{equation*}
where $V^{M(p)}$, for instance,
represents the generalized Verma module of level $\ell\in\N$ induced from the
$(p+1)$-dimensional irreducible $\mathfrak{sl}(2,\C)$-module $M(p)$ for $p \in \N$. The main difficulty in constructing an intertwining operator $\mathcal{Y}$ from an $\mathfrak{sl}(2,\C)$-homomorphism $f: M(p)\otimes M(q)\rightarrow M(r)$ is to first construct
\begin{equation*}
 \mathcal{Y}: M(p)\otimes M(q)\rightarrow V^{M(r)}\lbrace x\rbrace.
\end{equation*}
satisfying appropriate properties, where $V^{M(r)}\lbrace x\rbrace$ denotes the space of formal series involving general complex powers of $x$ with coefficients in $V^{M(r)}$. It is then comparatively routine to extend $\mathcal{Y}$ to an intertwining operator of type $\binom{V^{M(r)}}{V^{M(p)}\,V^{M(q)}}$.

Given a homomorphism $f$ from $M(p)\otimes M(q)$ into $M(r)$, there is a very
natural way to extend this to a map
$$\bar{\mathcal{Y}}: M(p)\otimes M(q)\rightarrow (V^{M(r)})'\lbrace x\rbrace $$
where $(V^{M(r)})'$ is the contragredient of $V^{M(r)}$. In fact, $\bar{\mathcal{Y}}$ extends to an intertwining operator of type $\binom{(V^{M(r)})'}{V^{M(p)}\,V^{M(q)}}$ by Li's results in \cite{Li3}. However, $(V^{M(r)})'\ncong V^{M(r)}$ since the contragredient of a reducible generalized Verma module is not a generalized Verma module. Nonetheless, there is a (degenerate) bilinear form on $V^{M(r)}$ whose radical is the maximal proper submodule $J(r)$. We use our theorem on the irreducibility of $J(r)$ to show that the image of $\bar{\mathcal{Y}}$ in $(V^{M(r)})'$ annihilates $J(r)$, so that we can use the bilinear form on $V^{M(r)}$ to transport $\bar{\mathcal{Y}}$ to a map
$$\mathcal{Y}^K: M(p)\otimes M(q)\rightarrow K(r)\lbrace x\rbrace$$
where $K(r)$ is an $\mathfrak{sl}(2,\C)$-module complement of $J(r)$ in $V^{M(r)}$. (Note that $K(r)$ is not an $\widehat{\mathfrak{sl}(2,\C)}$-module since $V^{M(r)}$ is indecomposable.)

Using an analogous procedure, again using the irreducibility of $J(r)$, we can
next construct a map
\begin{equation*}
 \mathcal{Y}^J: M(p)\otimes M(q)\rightarrow J(r)\lbrace x\rbrace
\end{equation*}
such that the operator
\[
\mathcal{Y}=\mathcal{Y}^K+\mathcal{Y}^J: M(p)\otimes M(q)\rightarrow V^{M(r)}\lbrace x\rbrace
\]
satisfies a commutator formula that then allows its extension to an intertwining operator of type $\binom{V^{M(r)}}{V^{M(p)}\,V^{M(q)}}$. In this way, we completely determine the space of intertwining operators of type
$\binom{V^{M(r)}}{V^{M(p)}\,V^{M(q)}}$. Additionally, we show that intertwining operators
of type $\binom{V^{M(r)}}{V^{M(p)}\,V^{M(q)}}$ induce intertwining operators of
type $\binom{L(r)}{V^{M(p)}\,L(q)}$ under certain easily-checked conditions.

Recalling the connection between vertex algebraic intertwining operators and
tensor products of modules for a vertex operator algebra, we expect that our
results here will contribute to solving the problem of constructing tensor
categories of $V^{M(\ell\Lambda_0)}$-modules. In particular, our determination
of the fusion rules can be interpreted as calculations of the
multiplicities of the modules $V^{M(r)}$ in tensor products of $V^{M(p)}$ and
$V^{M(q)}$. It would be interesting to see how our results here can be
generalized to larger classes of $\widehat{\mathfrak{sl}(2,\C)}$-modules, or to
generalized Verma modules for higher-rank affine Kac-Moody Lie algebras. The
primary restriction on our methods here is that they require knowing the
structure of generalized Verma modules. For instance, we cannot expect maximal
proper submodules of generalized Verma modules for higher-rank Kac-Moody Lie
algebras to be irreducible. Additional research is currently under way to
develop methods to bypass this problem.

This paper is structured as follows. In Section 2 we recall constructions and results that we shall need from the theory of affine Kac-Moody Lie algebras.  In Section 3, we recall the Kazhdan-Lusztig multiplicity formula and the Garland-Lepowsky resolutions that we shall then use to prove the irreducibility of maximal proper submodules for generalized Verma modules for
$\widehat{\mathfrak{sl}(2,\C)}$. In Section 4, we recall the vertex operator algebra and module structures on generalized Verma modules for an affine Lie algebra as well as the definition of intertwining operator among modules for a vertex operator algebra. In Section 5, we review results on bilinear pairings between generalized Verma modules that we shall need in Section 6, where we prove our main theorem on intertwining operators among generalized Verma modules for $\widehat{\mathfrak{sl}(2,\C)}$. In Section 7, we prove conditions
under which an intertwining operator of type $\binom{V^{M(r)}}{V^{M(p)}\,V^{M(q)}}$ descends to an intertwining operator of type $\binom{L(r)}{V^{M(p)}\,L(q)}$, and in Section 8 we construct examples of intertwining operators among $\widehat{\mathfrak{sl}(2,\C)}$-modules using our theorems. We also show by counterexample that the conditions on our theorems are necessary. Finally, in the Appendix, we give a proof that a map
\begin{equation*}
 \mathcal{Y}: M(p)\otimes M(q)\rightarrow V^{M(r)}\lbrace x\rbrace
\end{equation*}
satisfying appropriate conditions extends to an intertwining operator
\begin{equation*}
 \mathcal{Y}: V^{M(p)}\otimes V^{M(q)}\rightarrow V^{M(r)}\lbrace x\rbrace.
\end{equation*}
Note that although our main concern in this paper is with
$\widehat{\mathfrak{sl}(2,\C)}$, some results are proved for general untwisted
affine Kac-Moody algebras.

\paragraph{Acknowledgments}
We would like to express our deepest thanks to our advisor James Lepowsky for
inspiring us to think about this question. We are also grateful to Drazen Adamovic, Antun Milas
and Mirko Primc for inviting us to give talks at Representation Theory XIII
in Dubrovnik, Croatia in June 2013. We thank Antun Milas for informing us about Rocha-Caridi and Wallach's multiplicity formula for rank $2$ symmetrizable Kac-Moody Lie algebras in \cite{RW2} and Haisheng Li for explaining details
to us on his method of determining fusion rules among certain modules, and we
also thank Shashank Kanade for helpful discussions.

\section{Affine Kac-Moody Lie algebras}

In this section we shall recall material that we shall need on Kac-Moody Lie algebras. See references such as \cite{H}, \cite{Le1}, \cite{K}, and \cite{MP} for more details.

We take $\mathfrak{g}$ to be a finite-type Kac-Moody algebra with indecomposable Cartan matrix, that is, $\mathfrak{g}$ is a finite-dimensional simple Lie algebra over $\mathbb{C}$. If $\mathfrak{h}$ is a Cartan subalgebra of $\g$, then $\g$ has the triangular decomposition
\begin{equation*}
 \g=\g_-\oplus\mathfrak{h}\oplus\g_+,
\end{equation*}
where $\g_\pm$ are the sums of the positive and negative root spaces of $\g$, respectively. We use $\langle\cdot,\cdot\rangle $ to denote the symmetric nondegenerate invariant bilinear form  on $\mathfrak{g}$ such that long roots have square length $2$.

Using $\g$ and the bilinear form $\langle\cdot,\cdot\rangle$, we can construct the affine Lie algebra
\begin{equation*}
\widehat{\mathfrak{g}}=\mathfrak{g}\otimes\mathbb{C}[t,t^{-1}]\oplus\mathbb{C}
\mathbf{k}
\end{equation*}
with $\mathbf{k}$ central and all other bracket relations defined by
\begin{equation*}\label{affinecomm}
 [g\otimes t^m, h\otimes t^n]=[g,h]\otimes t^{m+n}+m\langle
g,h\rangle\delta_{m+n,0}\mathbf{k}
\end{equation*}
for $g,h\in\mathfrak{g}$ and $m,n\in\mathbb{Z}$. We have the decomposition
\begin{equation*}
\widehat{\mathfrak{g}}=\widehat{\mathfrak{g}}_-\oplus\widehat{\mathfrak{g}}
_0\oplus\widehat{\mathfrak{g}}_+
\end{equation*}
where $\widehat{\mathfrak{g}}_{\pm}=\coprod_{n\in\pm\mathbb{Z}_+} \mathfrak{g}\otimes t^n$ and $\widehat{\mathfrak{g}}_0=\mathfrak{g}\oplus\mathbb{C}\mathbf{k}$.

We now construct generalized Verma modules for $\ghat$, in the sense of \cite{GL} and \cite{Le2}. For a dominant integral weight $\lambda$ of $\mathfrak{g}$ and $\ell\in\mathbb{C}$, we take $M(\lambda,\ell)$ to be the $\widehat{\mathfrak{g}}_0$-module which is the
irreducible $\mathfrak{g}$-module with highest weight $\lambda$ on which
$\mathbf{k}$ acts as $\ell$. We then have the generalized Verma module
\begin{equation*}
V^{M(\lambda,\ell)}=U(\widehat{\mathfrak{g}})\otimes_{U(\widehat{\mathfrak{g}}
_0\oplus\widehat{\mathfrak{g}}_+)} M(\lambda,\ell)\cong
U(\widehat{\mathfrak{g}}_-)\otimes_\mathbb{C} M(\lambda,\ell),
\end{equation*}
where the linear isomorphism follows from the Poincar${\rm
\acute{e}}$-Birkhoff-Witt theorem.
 The scalar $\ell$ is called the \textit{level} of $V^{M(\lambda,\ell)}$. For $g\in\mathfrak{g}$ and $n\in\mathbb{Z}$, we use the notation $g(n)$ to denote the action of $g\otimes t^n$ on a $\widehat{\mathfrak{g}}$-module. Then the
generalized Verma module $V^{M(\lambda,\ell)}$ is linearly spanned by vectors of the form
\begin{equation*}
 g_1(-n_1)\cdots g_k(-n_k) u
\end{equation*}
where $g_i\in\mathfrak{g}$, $n_i>0$, and $u\in M(\lambda,\ell)$.

We will also need the semidirect product Lie algebra
\begin{equation*}
 \widetilde{\mathfrak{g}}=\widehat{\mathfrak{g}}\rtimes\mathbb{C}\mathbf{d}
\end{equation*}
where $[\mathbf{d},\mathbf{k}]=0$ and
\begin{equation*}
 [\mathbf{d},g\otimes t^n]=n(g\otimes t^n)
\end{equation*}
for $g\in\g$, $n\in\Z$. The Lie algebra $\widetilde{\mathfrak{g}}$ is the \textit{Kac-Moody Lie algebra} associated with a certain generalized Cartan matrix formed by adding one extra row and column to the Cartan matrix of $\mathfrak{g}$. The Kac-Moody algebra $\gtilde$ has the triangular decomposition
\begin{equation*}
 \gtilde=\gtilde_-\oplus\mathfrak{H}\oplus\gtilde_+
\end{equation*}
where $\gtilde_\pm=\g_\pm\oplus\ghat_\pm$ and
\begin{equation*}
 \mathfrak{H}=\mathfrak{h}\oplus\C\mathbf{k}\oplus\C\mathbf{d}
\end{equation*}
is a Cartan subalgebra of $\gtilde$. If $h_1,\ldots,h_l\in\mathfrak{h}$ are the coroots of $\g$, then $\mathfrak{H}$ has a basis $h_0, h_1,\ldots,h_l,\mathbf{d}$ where $h_0=-h_\theta+\mathbf{k}$ and $h_\theta$ denotes the coroot associated to the longest root $\theta$ of $\g$.

Note that
\begin{equation*}
 \mathfrak{H}^*=\mathfrak{h}^*\oplus\C\mathbf{k}'\oplus\C\mathbf{d}'
\end{equation*}
where $\mathbf{k}'$ and $\mathbf{d}'$ are dual to $\mathbf{k}$ and $\mathbf{d}$, respectively. The simple roots of $\gtilde$ are given by $\alpha_0,\alpha_1,\ldots,\alpha_l$, where $\alpha_1,\ldots,\alpha_l\in\mathfrak{h}^*$ are the simple roots of $\g$ and $\alpha_0=-\theta+\mathbf{d}'$. We recall the action of the Weyl group $\mathcal{W}$ of $\gtilde$ on $\mathfrak{H}^*$: $\mathcal{W}$ is the group generated by the simple reflections $r_0,r_1,\ldots,r_l$ where
\begin{equation*}
 r_i(\Lambda)=\Lambda-\Lambda(h_i)\alpha_i
\end{equation*}
for $\Lambda\in\mathfrak{H}^*$ and $0\leq i\leq l$.

Recall the set of integral weights of $\gtilde$
\[P = \{\Lambda \in \mathfrak{H}^*\, \vert\, \Lambda(h_i) \in \mathbb{Z}\ \
{\rm for}\  {\rm all}\  0\leq i\leq l\},\] and the set of dominant
integral weights
\[P^+ = \{\Lambda \in \mathfrak{H}^*\ |\
\Lambda(h_i) \in \mathbb{N}\ \  {\rm for}\  {\rm all}\   0\leq i\leq l\}.\]
We will in particular fix $\rho\in P^+$ such that $\rho(h_i)=1$ for $0\leq i\leq l$. We will also need the set
\begin{equation*}
 P_1=\lbrace \Lambda\in\mathfrak{H}^*\,\vert\,\Lambda(h_i)\in\N\ \
{\rm for}\  {\rm all}\  1\leq i\leq l\rbrace.
\end{equation*}
\begin{rema}
 Since the coroots of $\gtilde$ are contained in $\mathfrak{h}\oplus\C\mathbf{k}$, if $\Lambda\in P^+$, so is $\Lambda+h\mathbf{d}'$ for any $h\in\C$. Thus we can define the dominant integral weights of $\ghat$ to be the elements $\Lambda\in\mathfrak{h}^*\oplus\C\mathbf{k}'$ such that $\Lambda(h_i)\in\N$ for $0\leq i\leq l$. Then
 \begin{equation*}
  \Lambda=\lambda+\ell\mathbf{k}'
 \end{equation*}
is a dominant integral weight of $\ghat$ if and only if $\ell\in\N$ and $\lambda$ is a dominant integral weight of $\g$ satisfying $\langle\lambda,\theta\rangle\leq\ell$.
\end{rema}

For $\Lambda\in\mathfrak{H}^*$, let $\C_\Lambda$ denote the one-dimensional $\mathfrak{H}\oplus\gtilde_+$-module on which $\gtilde_+$ acts trivially and $\mathfrak{H}$ acts according to $\Lambda$. We define the \textit{Verma module}
\begin{equation*}
 V^\Lambda=U(\gtilde)\otimes_{U(\mathfrak{H}\oplus\gtilde_+)}\C_\Lambda\cong U(\gtilde_-)\otimes_\C \C_\Lambda.
\end{equation*}
We can also give generalized Verma modules for $\ghat$ a $\gtilde$-module structure. Note that for any $\Lambda\in P_1$ we can write
\begin{equation}\label{P1exam}
 \Lambda=\lambda+\ell\mathbf{k}'+h\mathbf{d}'
\end{equation}
where $\lambda$ is a dominant integral weight of $\g$ and $\ell,h\in\C$. Then $M(\lambda,\ell)$ becomes a $\ghat_0\oplus\C\mathbf{d}$-module on which $\mathbf{d}$ acts as $h$, and $V^{M(\lambda,\ell)}$ becomes a $\gtilde$-module via
\begin{equation*}
\mathbf{d}\cdot g_1(-n_1)\cdots g_k(-n_k) u=(h-n_1-\ldots-n_k) g_1(-n_1)\cdots g_k(-n_k) u
\end{equation*}
for $g_i\in\mathfrak{g}$, $n_i>0$, and $u\in M(\lambda,\ell)$. Note that as a $\gtilde$-module, $V^{M(\lambda,\ell)}$ is generated by a highest weight vector of weight $\Lambda$.
\begin{rema}
 We will typically use the notation $V^{M(\Lambda)}=V^{M(\lambda,\ell)}$ for $\Lambda$ as in \eqref{P1exam} when we need to consider a generalized Verma module as a $\gtilde$-module.
\end{rema}

For $\Lambda\in\mathfrak{H}^*$, the Verma module $V^\Lambda$ is the universal highest weight $\gtilde$-module with highest weight $\Lambda$; thus for $\Lambda\in P_1$, there is a surjection $\eta: V^\Lambda\rightarrow V^{M(\Lambda)}$ taking a basis vector $v_0$ of $\C_\Lambda$ to a highest weight vector generating $V^{M(\Lambda)}$. Similarly, for $\Lambda\in P_1$, $V^{M(\lambda)}$ is the universal highest weight $\gtilde$-module such that for any $\ghat_0\oplus\ghat_+\oplus\C\mathbf{d}$-module homomorphism
from $M(\Lambda)$ into a $\gtilde$-module $W$, there is a unique
extension to a $\gtilde$-module homomorphism from $V^{M(\Lambda)}$ into
$W$. Both $V^\Lambda$ and $V^{M(\Lambda)}$ have unique maximal proper submodules, the sum of all submodules that do not intersect the highest weight space. Since $V^{M(\Lambda)}$ is a quotient of $V^\Lambda$, both have the same unique irreducible quotient $L(\Lambda)$. We will use the notation $J(\Lambda)$ for the maximal proper submodule of $V^{M(\Lambda)}$.

\begin{rema}
 For $\Lambda\in P_1$ as in \eqref{P1exam}, we will also use the notation $J(\lambda)$ for $J(\Lambda)$ and $L(\lambda)$ for $L(\Lambda)$, especially if we are considering these modules as $\ghat$-modules and if $\ell$ is understood.
\end{rema}

Recall that a $\gtilde$-module $W$ is called {\it $\mathfrak{H}$-semisimple} if
\begin{equation*}
 W = \coprod_{\Lambda \in \mathfrak{H}^{*}}W_{\Lambda},
\end{equation*}
where
\begin{equation*}
 W_{\Lambda} = \{w \in W \ |\ h\cdot w = \Lambda(h)w\ {\rm for}\ h \in
\mathfrak{H}\}.
\end{equation*}
The subspace $W_{\Lambda}$ is called a {\it weight space}, and $\Lambda \in
\mathfrak{H}^{*}$ is called a {\it weight} if $W_{\Lambda} \neq 0$. We use
$Q^+$ to denote the subset of $\mathfrak{H}^*$ consisting of non-negative
integral combinations of the simple roots $\alpha_0,\ldots,\alpha_l$. Recall
that $\mathfrak{H}^*$ is partially ordered via $\alpha\leq\beta$ if and only if $\beta-\alpha\in Q^+$. For $\beta \in \mathfrak{H}^{*}$, set $D(\beta) =
\{\alpha \in \mathfrak{H}^*\ |\ \alpha \leq \beta\}$.

We recall the Bernstein-Gelfand-Gelfand category $\mathcal{O}$  whose objects
are $\gtilde$-modules $W$ satisfying:
\begin{enumerate}
\item[(i)] $W$ is $\mathfrak{H}$-semisimple with finite dimensional weight
spaces.
\item[(ii)] There exist finitely many elements $\beta_1, \dots \beta_k \in
\mathfrak{H}^*$ such that any weight $\Lambda$ of $M$ belongs to some
$D(\beta_i)$.
\end{enumerate}
The category $\mathcal{O}$ is stable under the operations of taking submodules, quotients, and finite direct sums. All highest weight $\gtilde$-modules, including Verma modules and generalized Verma modules, are in category $\mathcal{O}$. We shall need the following proposition:

\begin{prop}\label{p1} Suppose $W$ is a module in category $\mathcal{O}$ and
$\Lambda \in \mathfrak{H}^*$. Then there exists an increasing filtration $0
= W_0 \subseteq W_1 \subseteq \cdots \subseteq W_t = W$ of submodules of $W$ and
a subset $J$ of $\{1, \dots, t\}$ such that
\begin{enumerate}
\item[(i)]For $j \in J$, $W_j/W_{j-1} \cong L(\Lambda_j)$ for some $\Lambda_j
\geq \Lambda$ and
\item[(ii)]For $j \in \{1, \dots, t\}-J$ and any $\mu \geq \Lambda$,
$(W_j/W_{j-1})_{\mu} = 0$.
\end{enumerate}
\end{prop}

Such a filtration is called a {\it local composition series} of $W$ at $\Lambda$. Suppose $W$ is in category $\mathcal{O}$ and $\mu \in \mathfrak{H}^*$. Fix $\Lambda \in \mathfrak{H}^*$ such that $\mu \geq \Lambda$ and construct a local composition series of $W$ at $\Lambda$. Denote by $[W : L(\mu)]$ the number of times $\mu$ among $\{\Lambda_j \ |\ j\in J\}$. In fact $[W : L(\mu)]$ is independent of the filtration furnished by Proposition \ref{p1} and the choice of $\Lambda$. This number is called the {\it multiplicity} of $L(\mu)$ in $W$.

The main results in this paper apply to $\g=\mathfrak{sl}(2,\C)$, so we describe how some of the constructions in this section apply to this case and fix notation that we will use for $\mathfrak{sl}(2,\C)$. We take the usual basis $\lbrace h, e, f\rbrace$ of
$\mathfrak{sl}(2,\mathbb{C})$ satisfying the bracket relations
\[[h, e] = 2e,\ [h,f] = -2f,\ [e,f] = h,\]
and we take the symmetric invariant bilinear form on
$\mathfrak{sl}(2,\mathbb{C})$ to be given by
\[ \langle e,f\rangle=\dfrac{1}{2}\langle h,h\rangle =1,\ \langle e,e\rangle
=\langle f,f\rangle =\langle e,h\rangle=\langle f,h\rangle =0.\]
Then $\widetilde{\mathfrak{sl}(2,\mathbb{C})}$ has the Cartan subalgebra
\begin{equation*}
 \mathfrak{H}=\mathbb{C}h\oplus\mathbb{C}\mathbf{k}\oplus\mathbb{C}\mathbf{d},
\end{equation*}
and we have $h_0=-h+\mathbf{k}$ and $h_1=h$. We also consider
\begin{equation*}
\mathfrak{H}^*=\mathbb{C}\alpha\oplus\mathbb{C}\mathbf{k'}\oplus\mathbb{C}
\mathbf{d'},
\end{equation*}
where $\lbrace \alpha/2,\mathbf{k'},\mathbf{d'}\rbrace $ is the basis of
$\mathfrak{H}^*$ dual to $\lbrace h,\mathbf{k},\mathbf{d}\rbrace$ The simple roots of $\widehat{\mathfrak{sl}(2,\mathbb{C})}$ are
$\alpha_0=-\alpha+\mathbf{d'}$ and $\alpha_1=\alpha $. We may also take $\rho =\dfrac{\alpha}{2}+2\mathbf{k}'$ since then $\rho(h_0)=\rho(h_1)=1$. The Weyl group generators $r_0$ and $r_1$ act on $\mathfrak{H}^*$ via
\begin{align}\label{sl2Weyl}
r_0(\alpha)=\alpha-\alpha(-h+\mathbf{k})(-\alpha+\mathbf{d'})=-\alpha+2\mathbf{d'}
; &\hspace{2em} r_1(\alpha)=-\alpha;\nonumber\\
r_0(\mathbf{k'})=\mathbf{k'}-\mathbf{k'}(-h+\mathbf{k})(-\alpha+\mathbf{d'}
)=\alpha+\mathbf{k'}-\mathbf{d'}; &\hspace{2em}
r_1(\mathbf{k'})=\mathbf{k'};\nonumber\\
r_0(\mathbf{d'})=\mathbf{d'}; &\hspace{2em} r_1(\mathbf{d'})=\mathbf{d'}.
\end{align}

When $\g=\mathfrak{sl}(2,\C)$, the elements of $P_1$ have the form
\begin{equation*}
 \Lambda=n\dfrac{\alpha}{2}+\ell\mathbf{k}'+h\mathbf{d}'
\end{equation*}
where $n\in\N$. We will typically use $M(n,\ell)$, or $M(n)$ if $\ell$ is understood, to denote the finite-dimensional irreducible $\ghat_0\oplus\C\mathbf{d}$-module $M(\Lambda)$ of highest weight $\Lambda$. Thus we typically use $V^{M(n)}$ to refer to the generalized Verma module for $\widehat{\mathfrak{sl}(2,\C)}$ (or $\widetilde{\mathfrak{sl}(2,\C)}$) generated by $M(\Lambda)$. Note that as an $\mathfrak{sl}(2,\C)$-module, $M(n)$ is the irreducible $\mathfrak{sl}(2,\C)$-module of dimension $n+1$.

\setcounter{equation}{0}
\section{Irreducibility of maximal proper submodules in generalized Verma
modules for $\widehat{\mathfrak{sl}(2,\mathbb{C})}$}

In this section, we prove that for certain weights $\Lambda$ of
$\widehat{\mathfrak{sl}(2,\mathbb{C})}$, the maximal proper submodule of
$V^{M(\Lambda)}$ is irreducible. We use a formula of Rocha-Caridi and Wallach for multiplicities of irreducible modules in Verma modules as well as the resolutions of standard modules by generalized Verma modules of Garland and Lepowsky. 

We start by recalling the definition of a Coxeter system, examples of which are provided by the Weyl groups of Kac-Moody Lie algebras and the simple reflections which generate them:
\begin{defn}{\rm A {\it Coxeter system} is a pair $(\mathcal{W}, S)$ where $\mathcal{W}$ is a group and $S$ is a set of involutions in $\mathcal{W}$ such that $\mathcal{W}$ has a presentation of the form
\[
W = \langle S\;|\; (st)^{m(s, t)} \rangle.
\]
Here $m(s, t)$ denotes the order of the element $st$ in $\mathcal{W}$ and in the presentation for $\mathcal{W}$, $(s, t)$ ranges over all pairs in $S\times S$ such that $m(s, t) \neq \infty$. We further require the set $S$ to be finite. The group $\mathcal{W}$ is a {\it Coxeter group} and $S$ is a {\it fundamental set of generators} of $\mathcal{W}$.}
\end{defn}

\begin{defn}\label{d1}{\rm For a Coxeter system $(\mathcal{W},S)$, a \textit{reduced expression} for $\sigma\in\mathcal{W}$ in terms of the generators $S$ is an expression $\sigma=s_1\cdots s_t$, where $s_i\in S$, for which $t$ is minimal; we say $t$ is the \textit{length} of $\sigma$. The \textit{Bruhat order} on $\mathcal{W}$ is defined by $\sigma\leq\tau$ for $\sigma,\tau\in \mathcal{W}$ if a reduced expression for $\sigma$ in terms of the generators $S$ occurs as a subword of a reduced expression for $\tau$. }
\end{defn}

The following multiplicity result for $\widetilde{\mathfrak{sl}(2, \mathbb{C})}$ is Corollary 4.14(i) in \cite{RW2}, which applies more generally to any rank $2$ symmetrizable Kac-Moody Lie algebra:
\begin{thm}\label{c1}
For  $\Lambda$ a dominant integral weight of $\widetilde{\mathfrak{sl}(2,
\mathbb{C})}$ and $x, y$ in the Weyl group of $\widetilde{\mathfrak{sl}(2,
\mathbb{C})}$,
 \[[V^{x(\Lambda + \rho) - \rho}: L(y(\Lambda + \rho) - \rho)] = \left\{
\begin{array}{ll} 1\ \ {\rm if}\ \ x \leq y;\\ 0 \ \ {\rm otherwise}.
\end{array} \right.\]
\end{thm}

\begin{rema}
 Theorem \ref{c1} can be generalized to any symmetrizable Kac-Moody Lie algebra $\mathfrak{G}$ using Kazhdan-Lusztig polynomials. In \cite{KL}, Kazhdan and Lusztig constructed polynomials $P_{x,y}(q)\in\Z[q]$, where $x$ and $y$ are elements of some Coxeter group $\mathcal{W}$, using the Hecke algebra associated to $\mathcal{W}$. If $\mathcal{W}$ is the Weyl group of $\mathfrak{G}$, Deodhar, Gabber, and Kac conjectured (\cite{DGK}) that
 \begin{equation}\label{KLmult}
  [V^{x(\Lambda + \rho) - \rho}: L(y(\Lambda + \rho) - \rho)] = P_{x,
y}(1)
 \end{equation}
for any $\Lambda \in P^{+}$ and $x, y\in\mathcal{W}$. This conjecture was proved by Kashiwara-Tanisaki (\cite{Ka}, \cite{KaT}) and Casian (\cite{C}). In general, $P_{x,y}(q)=0$ unless $x\leq y$, and if the Coxeter group $\mathcal{W}$ has rank $2$, $P_{x,y}(q)=1$ for $x\leq y$. Thus in the case that $\mathfrak{G}$ has rank $2$, \eqref{KLmult} reduces to Theorem \ref{c1}. 
\end{rema}

To obtain multiplicity results for generalized Verma modules for $\widetilde{\mathfrak{sl}(2,\C)}$, we need the following analogue for $\widehat{\mathfrak{sl}(2, \C)}$ of Proposition 2.1 in \cite{Le2}:
\begin{prop}\label{cc}
Let $\Lambda \in P_1$ and let $v_0$ be a highest weight vector of the Verma module $V^{\Lambda}$. Then there is an $\widehat{\mathfrak{sl}(2, \C)}$-module exact sequence
\[
0 \longrightarrow V^{r_1(\Lambda + \rho) - \rho}
\stackrel{\xi}{\longrightarrow} V^{\Lambda}
\stackrel{\eta}{\longrightarrow} V^{M(\Lambda)} \longrightarrow
0,
\]
where $\eta$ takes $v_0$ to a highest weight vector generating
$V^{M(\Lambda)}$, and $\xi$ takes a highest weight vector generating
$V^{r_1(\Lambda + \rho) - \rho}$ to a nonzero multiple of the
highest weight vector $f^{\Lambda(h) + 1}\cdot v_0 \in
V^{\Lambda}$.
\end{prop}
\pf The fact that $\xi$ is injective follows from \cite{RW1} (cf. \cite{V}).
We only need to prove that ${\rm Ker}\; \eta$ is the $\ghat=\widehat{\mathfrak{sl}(2, \mathbb{C})}$-submodule of $V^{\lambda}$ generated by
$f^{\Lambda(h) + 1}\cdot v_{\Lambda}$; we denote this module by
$K^{\Lambda}$.

Let and let $X$ be the $\ghat_0$-submodule of $V^{\Lambda}$ generated by $v_{\Lambda}$. Then by a classical theorem of Harish-Chandra (see for instance Section 21.4 of \cite{H}), the kernel of the surjection $\eta|_X: X \rightarrow M(\Lambda)$ is the $\ghat_0$-submodule of $X$ generated by $f^{\lambda(h) + 1}\cdot v_{\Lambda}$; here we regard $M(\Lambda)$ in the natural way as a $\ghat_0$-submodule of $V^{M(\Lambda)}$. Thus $K^{\Lambda} \subseteq {\rm
Ker}\; \eta$.

On the other hand, there are $\ghat_0$-module maps $M(\Lambda) \rightarrow
X/{\rm Ker}\; (\eta|_X) \rightarrow V^{\Lambda}/K^{\Lambda}$ whose composition
is nonzero, and so we have a natural $\ghat_0$-module injection $\iota:
M(\Lambda) \rightarrow V^{\Lambda}/K^{\Lambda}$. The image of this injection is clearly annihilated by $\ghat_+$, so that $\iota$ extends to a
$\ghat=\widehat{\mathfrak{sl}(2,\C)}$-module map $V^{M(\Lambda)} \rightarrow V^{\Lambda}/K^{\Lambda}$, by the universal property of generalized Verma modules. Since $V^{M(\Lambda)} = V^{\Lambda}/{\rm Ker}\;\eta$, we see that ${\rm Ker}\;\eta \subseteq K^{\Lambda}$.\epfv

We now combine Theorem \ref{c1} and Proposition \ref{cc} to obtain:
\begin{cor}\label{multiplicity 0}
Let $\gtilde = \widetilde{\mathfrak{sl}(2, \mathbb{C})}$, and let $\Lambda
\in P^{+}$, $x, y \in \mathcal{W}$ be such that $x(\Lambda + \rho) -~ \rho,\linebreak
y(\Lambda +~ \rho) - \rho \in P_1$. Then
\[[V^{M(x(\Lambda + \rho) - \rho)}: L(y(\Lambda + \rho) - \rho)] = \left\{
\begin{array}{ll} 1\ \ {\rm if}\ \ x \leq y,\ r_1x \nleq y;\\ 0 \ \ {\rm
otherwise}.
\end{array} \right.\]
\end{cor}
\pf For the case $x \nleq y$, multiplicity $0$ follows easily from Theorem
\ref{c1}. Thus suppose $x \leq y$, so that $[V^{x(\Lambda + \rho) - \rho}:
L(y(\Lambda + \rho) - \rho)]=1$. Since by Proposition \ref{cc} $V^{M(x(\Lambda +
\rho) - \rho)} \cong V^{x(\Lambda + \rho) - \rho}/V^{r_1x(\Lambda + \rho) -
\rho}$, we can get a local composition of $V^{x(\Lambda + \rho) - \rho}$ at
$y(\Lambda + \rho) - \rho$ by joining together local composition series of
$V^{r_1x(\Lambda + \rho) - \rho}$ and  $V^{M(x(\Lambda + \rho) - \rho)}$ at
$y(\Lambda+\rho)-\rho$. If $r_1x\leq y$, then by Theorem \ref{c1}, $
L(y(\Lambda + \rho) - \rho)$ appears once in any local composition series for
$V^{r_1x(\Lambda + \rho) - \rho}$, and so it cannot appear in a local
composition series for $V^{M(x(\Lambda + \rho) - \rho)}$. On the other hand, if
$r_1x\nleq y$, then by Theorem \ref{c1} $ L(y(\Lambda + \rho) - \rho)$ does
not appear in any local composition series for $V^{r_1x(\Lambda + \rho) -
\rho}$, and so it must appear once in any local composition series for
$V^{M(x(\Lambda + \rho) - \rho)}$. \epfv

Suppose that $\Lambda$ is a dominant integral weight of
$\widehat{\mathfrak{sl}(2,\mathbb{C})}$. We now recall the
$\widehat{\mathfrak{sl}(2, \mathbb{C})}$-module resolution for the standard
module $L(\Lambda)$ from \cite{GL}. Let $\mathcal{W}^1$ be the subset of the Weyl group $\mathcal{W}$ consisting of those $w
\in \mathcal{W}$ such that $w^{-1}\alpha_1$ is a positive root. Then we can easily see that
\[\mathcal{W}^1 = \{(r_0r_1)^nr_0, (r_0r_1)^n\ |\ n \in \N\}.\]
Let $\mathcal{W}_j^1$ be the set of elements of $\mathcal{W}^1$ of length $j$, for $j \in\N$. Then the set $\mathcal{W}_j^1$ only consists of one element $w_j$, that is,
\begin{eqnarray*}
w_j = \left\{
\begin{array}{lll}(r_0r_1)^nr_0 & \mathrm{if} & j = 2n+1
\\ (r_0r_1)^n & \mathrm{if} & j = 2n
\end{array} \right. .
\end{eqnarray*}
The next theorem gives an $\widehat{\mathfrak{sl}(2, \C)}$-module resolution for
$L(\Lambda)$ which follows from Theorem $8.7$ in (\cite{GL}):

\begin{thm}\label{resolution} There is an $\widehat{\mathfrak{sl}(2,
\C)}$-module
resolution for $L(\Lambda)$ as follows:
\[\cdots \rightarrow V^{M(w_j(\Lambda + \rho) - \rho)}
\stackrel{d_j}{\rightarrow}
V^{M(w_{j-1}(\Lambda + \rho) -\rho)}
\stackrel{d_{j-1}}{\rightarrow} \cdots
\stackrel{d_2}{\rightarrow} V^{M(r_0(\Lambda +\rho) - \rho)}
\stackrel{d_1}{\rightarrow} V^{M(\Lambda)}
\stackrel{\Pi_{\Lambda}}{\rightarrow} L(\Lambda) \rightarrow
0,\] where $d_j$ is a $\widehat{\mathfrak{sl}(2, \C)}$-module map from
$V^{M(w_j(\Lambda + \rho) - \rho)}$ to $V^{M(w_{j-1}(\Lambda + \rho) - \rho)}$,
for $j\geq 0$.\end{thm}

Now we prove the main theorem of this section:
\begin{thm}\label{irred}
In the setting of Theorem \ref{resolution},  $\mathrm{Im}\,d_j$ is an
irreducible and maximal proper submodule of $V^{M(w_{j-1}(\Lambda + \rho) - \rho)}$.
More precisely, we have
\[
{\rm Im}\ d_j\; =\; J(w_{j-1}(\Lambda + \rho) - \rho)\; \cong \; L(w_{j}(\Lambda
+ \rho) - \rho),
\]
for $j \in \Z_{+}$.
\end{thm}
\pf It follows from Theorem 2 in \cite{KK} (cf. Corollary 4.13 in \cite{RW2}) that if $\Lambda$ is a dominant integral
weight of $\widehat{\mathfrak{sl}(2,\mathbb{C})}$ and $L(\mu)$ appears in a
local composition series for $V^{M(w(\Lambda+\rho)-\rho)}$ where $w\in W^1$,
then $\mu=w'(\Lambda+\rho)-\rho\in P_1$ for some $w'\in W$. Then by Corollary \ref{multiplicity 0}, the generalized
Verma module $V^{M(w_{j-1}(\Lambda + \rho) - \rho)}$ only contains highest
weight vectors of weight $w_{j-1}(\Lambda + \rho) - \rho$ and weight
$w_{j}(\Lambda + \rho) - \rho$. Since Im $d_j$ is generated by a highest weight
vector of weight $w_{j}(\Lambda + \rho) - \rho$ in $V^{M(w_{j-1}(\Lambda + \rho)
- \rho)}$, any non-zero proper submodule of Im $d_j$ would provide a highest weight
vector of weight other than $w_{j-1}(\Lambda + \rho) - \rho$ or $w_{j}(\Lambda +
\rho) - \rho$ in $V^{M(w_{j-1}(\Lambda + \rho) - \rho)}$. Hence $\mathrm{Im}\,d_j$ has no
nontrivial submodule and is irreducible for any $j\geq 1$. Also, notice that since
\begin{equation*}
V^{M(w_{j-1}(\Lambda + \rho) - \rho)}/{\rm Im}\ d_{j} = V^{M(w_{j-1}(\Lambda + \rho) - \rho)}/{\rm Ker}\ d_{j-1}\cong \mathrm{Im}\,d_{j-1}
\end{equation*}
is irreducible, Im $d_j$ has to be maximal in $V^{M(w_{j-1}(\Lambda + \rho) - \rho)}$.
\epfv

\begin{rema}{\rm From Theorem \ref{irred}, we obtain the following exact
sequence
\[
0 \longrightarrow L(w_{j+1}(\Lambda + \rho) - \rho) \longrightarrow
V^{M(w_j(\Lambda + \rho) - \rho)}
\longrightarrow L(w_{j}(\Lambda + \rho) - \rho) \longrightarrow 0\]
for $j \in \N$.}
\end{rema}

\setcounter{equation}{0}
\section{Vertex operator algebras and modules from affine Lie algebras}

In this section $\g$ is any finite-dimensional simple Lie algebra over $\C$.
We recall the vertex algebraic structure on generalized Verma
modules for $\ghat$; see for example \cite{LL} for more details.

We use the definition of a vertex operator algebra $(V=\coprod_{n\in\mathbb{Z}} V_{(n)}, Y,\mathbf{1},\omega)$ from \cite{LL}, where the grading on $V$ is by conformal weight. We also use the definition of a module $(W=\coprod_{h\in\mathbb{C}} W_{(n)}, Y_W)$ for a vertex operator algebra $V$ from \cite{LL}. Note that we will typically drop the subscript from the module vertex operator $Y_W$ if the module $W$ is understood. If $V$ is any vertex operator algebra and $W$ is any $V$-module, the
\textit{contragredient} of $W$ is the $V$-module
\begin{equation*}
 W'=\coprod_{h\in\mathbb{C}} W_{(h)}^*
\end{equation*}
with vertex operator map given by
\begin{equation*}
 \langle Y_{W'}(v,x)w',w\rangle =\langle w',
Y_W(e^{xL(1)}(-x^{-2})^{L(0)}v,x^{-1})w\rangle,
\end{equation*}
where $\langle\cdot,\cdot\rangle$ denotes the pairing between a vector space and its dual. See \cite{FHL} for the proof that this gives a $V$-module structure on $W'$.

We fix a level $\ell\in\C$. When $\ell\neq -h^\vee$, where $h^\vee$ is the dual Coxeter number of $\g$, the generalized Verma module $V^{M(0,\ell)}$ (induced from the one-dimensional $\mathfrak{g}$-module $\mathbb{C}\mathbf{1}$), has the structure of a vertex operator algebra with vacuum $\mathbf{1}$ and vertex operator map determined by
\begin{equation}\label{vrtxop}
 Y(g(-1)\mathbf{1},x)=g(x)=\sum_{n\in\mathbb{Z}} g(n) x^{-n-1}
\end{equation}
for $g\in\mathfrak{g}$. The conformal vector $\omega$ of $V^{M(0,\ell)}$ is
given by
\begin{equation*}
 \omega=\dfrac{1}{2(\ell+h^\vee)}\sum_{i=1}^{\mathrm{dim}\,\mathfrak{g}} u_i(-1)
u_i(-1)\mathbf{1},
\end{equation*}
where $\lbrace u_i\rbrace$ is an orthonormal basis of $\mathfrak{g}$ with
respect to the form $\langle\cdot,\cdot\rangle$.

For any dominant integral weight $\lambda$ of $\mathfrak{g}$, the generalized Verma module
$V^{M(\lambda,\ell)}$ is a $V^{M(0,\ell)}$-module with vertex operator map also
determined by (\ref{vrtxop}). The conformal weight grading on
$V^{M(\lambda,\ell)}$ is given by
\begin{equation*}
 \mathrm{wt}\,g_1(-n_1)\cdots g_k(-n_k)u=n_1+\ldots
+n_k+\frac{\langle\lambda,\lambda+2\rho_{\mathfrak{g}}\rangle}{2(\ell+h^\vee)}
\end{equation*}
for $g_i\in\mathfrak{g}$, $n>0$, and $u\in M(\lambda,\ell)$, where
$\rho_{\mathfrak{g}}$ is half the sum of the positive roots of $\mathfrak{g}$.
We observe that $g(n)$ decreases weight by $n$, and so
\begin{equation}\label{l0withgn}
 [L(0), g(n)]=-n g(n)
\end{equation}
for any $g\in\mathfrak{g}$ and $n\in\mathbb{Z}$. In particular, $g(0)$ preserves
weights, so each weight space of $V^{M(\lambda,\ell)}$ is a (finite-dimensional)
$\mathfrak{g}$-module.

\begin{rema}
 Note that \eqref{l0withgn} implies that $V^{M(\lambda,\ell)}$ has a natural $\gtilde$-module structure on which $\mathbf{d}$ acts as $-L(0)$.
\end{rema}

\begin{rema}\label{gradingnotation}
 For convenience, we will shift the grading of any graded subspace $X$ of a
$V^{M(0,\ell)}$-module $W$ as follows: if the lowest conformal weight of $W$ is
some $h\in\mathbb{C}$, then we define $X(n)=X\cap W_{(n+h)}$, so that $X$ is
$\mathbb{N}$-graded:
 \begin{equation*}
  X=\coprod_{n\geq 0} X(n).
 \end{equation*}
\end{rema}

We remark that any $\widehat{\mathfrak{g}}$-submodule of $V^{M(\lambda)}$ is a $V^{M(0)}$-module as
well, and vice versa. Thus any quotient of $V^{M(\lambda)}$ is a
$V^{M(0)}$-module. In particular, the maximal proper submodule $J(\lambda)$ is a
$V^{M(0)}$-submodule, and the
irreducible quotient $L(\lambda)$ is a $V^{M(0)}$-module. Since the weight spaces of
$V^{M(\lambda)}$ are finite-dimensional $\mathfrak{g}$-modules and since
finite-dimensional $\mathfrak{g}$-modules are completely reducible, there is a
graded $\mathfrak{g}$-module $K(\lambda)$ such that
$V^{M(\lambda)}=K(\lambda)\oplus J(\lambda)$ (but $K(\lambda)$ certainly need
not be a $\widehat{\mathfrak{g}}$-module). We use $P_K$ and $P_J$ to refer to
the
projections from $V^{M(\lambda)}$ to $K(\lambda)$ and $J(\lambda)$,
respectively;
these are $\mathfrak{g}$-module homomorphisms, but not generally
$\widehat{\mathfrak{g}}$-module homomorphisms. There is a $\mathfrak{g}$-module
isomorphism between $K(\lambda)$ and $L(\lambda)$ which sends $u\in K(\lambda)$
to $u+J(\lambda)$.

The contragredient of a $V^{M(0)}$-module $W$ is also a
$\widehat{\mathfrak{g}}$-module, with the action of $\widehat{\mathfrak{g}}$
given by
\begin{equation*}
 \langle g(n)w',w\rangle=-\langle w',g(-n)w\rangle
\end{equation*}
for $g\in\mathfrak{g}$, $n\in\mathbb{Z}$, $w'\in W'$, and $w\in W$. Note that
this is not the $\widehat{\mathfrak{g}}$-module structure on $W'$ viewed as a
subspace of the dual module $W^*$. The module $W'$ is in category $\mathcal{O}$
(see for example \cite{MP}).

Since the goal of this paper is to construct intertwining operators among $V^{M(0)}$-modules, we recall the definition of intertwining operator among a triple of modules for a vertex operator algebra from \cite{FHL}.
For a general vector space $W$, we use $W\{x\}$ to denote the vector space of formal series of the form $\sum_{n \in \C} w_nx^n$, $w_n \in W$.
\begin{defn}{\rm
Let $W_1$, $W_2$ and $W_3$
be modules for a vertex operator algebra $V$. An {\it intertwining
operator} of type $\binom{W_3}{W_1\,W_2}$ is a linear map
\begin{eqnarray*}
\mathcal{Y}: W_1\otimes W_2&\to& W_3\{x\},
\\
w_{(1)}\otimes w_{(2)}&\mapsto &\mathcal{Y}(w_{(1)},x)w_{(2)}=\sum_{n\in
{\mathbb C}}(w_{(1)})_n
w_{(2)}x^{-n-1}\in W_3\{x\}
\end{eqnarray*}
satisfying the
following conditions:

\begin{enumerate}

\item  {\it Lower truncation}: For any $w_{(1)}\in W_1$, $w_{(2)}\in W_2$ and
$n\in
\mathbb{C}$,
\begin{equation*}
(w_{(1)})_{n+m}w_{(2)}=0\;\;\mbox{ for }\;m\in {\mathbb
N} \;\mbox{ sufficiently large.}
\end{equation*}

\item The {\it Jacobi identity}:
\begin{align*}
 x^{-1}_0\delta \bigg( \frac{x_1-x_2}{x_0}\bigg)
Y_{W_{3}}(v,x_1)\mathcal{Y}(w_{(1)},x_2)w_{(2)} & -  x^{-1}_0\delta \bigg(
\frac{-x_2+x_1}{x_0}\bigg)
\mathcal{Y}(w_{(1)},x_2)Y_{W_{2}}(v,x_1)w_{(2)}\nn
& =  x^{-1}_2\delta \bigg( \frac{x_1-x_0}{x_2}\bigg) \mathcal{
Y}(Y_{W_{1}}(v,x_0)w_{(1)},x_2) w_{(2)}
\end{align*}
for $v\in V$, $w_{(1)}\in W_1$ and $w_{(2)}\in W_2$.
\item The {\em $L(-1)$-derivative property:} for any
$w_{(1)}\in W_1$,
\begin{equation*}
\mathcal{Y}(L(-1)w_{(1)},x)=\frac d{dx}\mathcal{Y}(w_{(1)},x).
\end{equation*}
\end{enumerate}}
\end{defn}
\begin{rema}
We denote the vector space of intertwining operators of type
$\binom{W_3}{W_1\,W_2}$ by $\mathcal{V}^{W_3}_{W_1 W_2}$ and the corresponding
\textit{fusion rule} is given by $\mathcal{N}^{W_3}_{W_1
W_2}=\mathrm{dim}\,\mathcal{V}^{W_3}_{W_1 W_2}$.
\end{rema}
\begin{rema}\label{commit}
 Taking the coefficient of $x_0^{-1}$ in the Jacobi identity for intertwining
operators yields the \textit{commutator formula}
 \begin{equation*}
Y_{W_3}(v,x_1)\mathcal{Y}(w_{(1)},x_2)-\mathcal{Y}(w_{(1)},x_2)Y_{W_2}(v,
x_1)=\mathrm{Res}_{x_0}\,x^{-1}_2\delta \bigg( \frac{x_1-x_0}{x_2}\bigg)
\mathcal{
Y}(Y_{W_{1}}(v,x_0)w_{(1)},x_2)
 \end{equation*}
for $v\in V$ and $w_{(1)}\in W_1$. Similarly, taking the coefficient of
$x_1^{-1}$ and then the coefficient of $x_0^{-n-1}$ for $n\in\mathbb{Z}$ yields
the \textit{iterate formula}
\begin{equation*}
 \mathcal{Y}(v_n w_{(1)},x_2)=\mathrm{Res}_{x_1}\,\left((x_1-x_2)^n
Y_{W_3}(v,x_1)\mathcal{Y}(w_{(1)},x_2)-(-x_2+x_1)^n \mathcal{Y}(w_{(1)},x_2)
Y_{W_2}(v,x_1)\right)
\end{equation*}
for $v\in V$ and $w_{(1)}\in W_1$. The commutator and iterate formulas for $v$
and $w_{(1)}$ together are equivalent to the Jacobi identity for $v$ and
$w_{(1)}$ (see \cite{LL} for the special case where $\mathcal{Y}$ is a module
vertex operator).
\end{rema}

\begin{rema}
We recall from \cite{FHL} that if $W_1$, $W_2$, and $W_3$ are three
$V$-modules with lowest conformal weights contained in $h_1+\N$, $h_2+\N$, and $h_3+\N$,
respectively, then any $\mathcal{Y}\in \mathcal{V}^{W_3}_{W_1 W_2}$ can be written
\begin{equation}\label{intwop}
 \mathcal{Y}(u,x)v=\sum_{m\in\mathbb{Z}} \mathcal{Y}_m(u)v\,x^{-m+h_3-h_1-h_2},
\end{equation}
where for $u\in W_1(k)=(W_1)_{(h_1+k)}$ and $v\in W_2(n)=(W_2)_{(h_2+n)}$,  $\mathcal{Y}_m(u)v\in W_3(k+n-m)$. In particular, $\mathcal{Y}_0$ maps $W_1(0)\otimes W_2(0)$ into $W_3(0)$.
\end{rema}

\setcounter{equation}{0}
\section{Invariant pairings between generalized Verma modules}

In this section we continue to allow $\mathfrak{g}$ to be any finite-dimensional simple Lie
algebra over $\mathbb{C}$, and we work with $\widehat{\mathfrak{g}}$-modules of
a fixed level $\ell$. One of the main tools we will use to prove our main theorems on intertwining operators among $\widehat{\mathfrak{sl}(2,\C)}$-modules is a bilinear pairing  between
generalized Verma modules $V^{M(\lambda)}$ and $V^{M(\lambda^*)}$, where as
$\mathfrak{g}$-modules $M(\lambda^*)\cong M(\lambda)^*$. More precisely, we want
a bilinear map
\begin{equation*}
 \langle\cdot,\cdot\rangle_M : V^{M(\lambda^*)}\times
V^{M(\lambda)}\rightarrow\mathbb{C}
\end{equation*}
that is \textit{invariant} in the sense that
\begin{equation}\label{invariance}
 \langle g(m)u,v\rangle_M =-\langle u, g(-m)v\rangle_M
\end{equation}
for $g\in\mathfrak{g}$, $m\in\mathbb{Z}$, $u\in V^{M(\lambda^*)}$, and $v\in
V^{M(\lambda)}$. Note that we have a (nondegenerate) such pairing between
$(V^{M(\lambda)})'$ and $V^{M(\lambda)}$, but $(V^{M(\lambda)})'$ is not
isomorphic to $V^{M(\lambda^*)}$.
\begin{rema}
An invariant bilinear pairing between the generalized Verma modules
$V^{M(\lambda)}$ and $V^{M(\lambda^*)}$ can be induced from the Shapovalov
pairing between the corresponding Verma modules (see for example \cite{MP} for
information on Shapovalov forms). For convenience, we shall give direct proofs
of all the results that we need here, in the context of generalized Verma
modules.
\end{rema}
\begin{rema}
Some of the results in this section when $\lambda=0$ have also been proved in a
vertex algebraic setting in \cite{Li1}.
\end{rema}
\begin{rema}
 We will use the notation $\langle\cdot,\cdot\rangle_M$ to denote
the invariant bilinear pairing between generalized Verma modules
$V^{M(\lambda^*)}$
and $V^{M(\lambda)}$, and we will now reserve the notation
$\langle\cdot,\cdot\rangle$
to denote the pairing between a module and its contragredient, or the form on
$\mathfrak{g}$.
\end{rema}

To start with, we observe that there is an involution (anti-automorphism) of
$\widehat{\mathfrak{g}}$ determined by $g(m)\mapsto -g(-m)$ and
$\mathbf{k}\mapsto\mathbf{k}$. This extends to an involution $\theta$ of
$U(\widehat{\mathfrak{g}})$ that interchanges $U(\widehat{\mathfrak{g}}_+)$ and
$U(\widehat{\mathfrak{g}}_-)$. We can then define a bilinear pairing between
$V^{M(\lambda)}\cong U(\widehat{\mathfrak{g}}_-)\otimes M(\lambda)$ and
$V^{M(\lambda^*)}$ via
\begin{equation*}
 \langle u, y\otimes v\rangle_M =\langle P(\theta(y)u),v\rangle
\end{equation*}
for $u\in V^{M(\lambda^*)}$, $y\in U(\widehat{\mathfrak{g}}_-)$, and $v\in
M(\lambda)$, where $P$ denotes projection to the lowest conformal weight space
$M(\lambda^*)$ and $\langle\cdot,\cdot\rangle$ is any nondegenerate $\g$-invariant pairing between $M(\lambda^*)\cong M(\lambda)^*$ and $M(\lambda)$.
\begin{prop}\label{invariantform}
The bilinear form $\langle\cdot,\cdot\rangle_M$ is invariant.
\end{prop}
\pf
 Suppose $u\in V^{M(\lambda^*)}$, $y\in U(\widehat{\mathfrak{g}}_-)$, and $v\in
M(\lambda)$. For $g\in\mathfrak{g}$ and $m>0$, we have
\begin{align*}
 \langle g(m)u,y\otimes v\rangle_M & =  \langle\theta(y)\theta(-g(-m))
u,v\rangle_M \nonumber\\
& =  \langle \theta(-g(-m)y)u,v\rangle_M \nonumber\\
& = \langle u,-g(-m)y\otimes v\rangle_M,
\end{align*}
as desired. For $m=0$, we observe that for any $y\in
U(\widehat{\mathfrak{g}}_\pm)$, $[g(0),y]\in U(\widehat{\mathfrak{g}}_\pm)$ as
well. Thus, using the $\mathfrak{g}$-module actions on $M(\lambda)$ and
$M(\lambda^*)$, and the fact that $g(0)$ commutes with $P$ (since $g(0)$
preserves conformal weight), we obtain:
\begin{align*}
 \langle g(0)u,y\otimes
v\rangle_M&=\langle\theta(y)\theta(-g(0))u,v\rangle_M\nonumber\\
& =  \langle \theta(-g(0))\theta(y)u,v\rangle_M+\langle
[\theta(y),\theta(-g(0))]u,v\rangle_M\nonumber\\
& =  \langle
P(g(0)\theta(y)u),v\rangle+\langle\theta([-g(0),y])u,v\rangle_M\nonumber\\
& =  \langle g\cdot P(\theta(y)u),v\rangle+\langle u,[-g(0),y]\otimes
v\rangle_M\nonumber\\
& =  \langle P(\theta(y)u), -g\cdot v\rangle+\langle u,-g(0)y\otimes
v\rangle_M+\langle u,y\otimes g\cdot v\rangle_M\nonumber\\
& =  \langle u, -g(0)y\otimes v\rangle .
\end{align*}
Now, for $m<0$, we will prove invariance by induction on the degree of $y$ in
$y\otimes
v\in V^{M(\lambda)}$. For the base case $y=1$, we have
\begin{equation*}
 \langle g(m)u,v\rangle_M=0=\langle u,-g(-m)v\rangle_M
\end{equation*}
since $P(g(m)u)=0$ and $g(-m)v=0$. Now suppose invariance holds for any
$y\otimes v$ with the degree of $y$ less than some $k$. It is enough to show
that
\begin{equation*}
 \langle g(m) u,h(-n)y\otimes v\rangle_M=\langle  u,-g(-m)h(-n)y\otimes
v\rangle_M
\end{equation*}
when $h\in\mathfrak{g}$, $n>0$, and $n+\mathrm{deg}\,y=k$. Thus:
\begin{align*}
 \langle g(m) u,h(-n)y & \otimes v\rangle_M  =   -\langle h(n)g(m)u,y\otimes
v\rangle_M\nonumber\\
 & =  -\langle g(m)h(n)u,y\otimes v\rangle_M-\langle ([h,g](n+m)+n\langle
h,g\rangle\delta_{n+m,0}\ell)u,y\otimes v\rangle_M\nonumber\\
 & = \langle h(n)u,g(-m)y\otimes v\rangle_M+\langle u,([h,g](-n-m)-n\langle
h,g\rangle\delta_{n+m,0}\ell)y\otimes v\rangle_M\nonumber\\
 & =  -\langle u,h(-n)g(-m)y\otimes v\rangle_M+\langle u,
[h(-n),g(-m)]y\otimes
v\rangle_M\nonumber\\
 & =  \langle u,-g(-m)h(-n)y\otimes v\rangle_M,
\end{align*}
completing the proof of invariance.
\epfv
\begin{prop}\label{orthogonal}
 If $u\in V^{M(\lambda^*)}(m)$, $v\in V^{M(\lambda)}(n)$ and $m\neq n$, then
$\langle
u,v\rangle_M=0$.
\end{prop}
\pf
 If $u=y\otimes u'$ and $v=y'\otimes v'$ where $u'\in M(\lambda^*)$, $v'\in
M(\lambda)$, deg $y=m$ and deg $y'=n$, then
 \begin{equation*}
  \langle u,v\rangle_M=\langle y\otimes u',y'\otimes v'\rangle_M=\langle
P(\theta(y')y\otimes u'),v'\rangle_M.
 \end{equation*}
This is non-zero only if $\theta(y')y\otimes u'\in V^{M(\lambda^*)}(0)$, which
only
happens when deg $y'=n=m=$ deg $y$.
\epfv
\begin{prop}\label{radicals}
 The left and right radicals of $\langle\cdot,\cdot\rangle_M$ are the maximal
proper submodules $J(\lambda^*)$ and $J(\lambda)$, respectively.
\end{prop}
\pf
 Since $\langle\cdot,\cdot\rangle_M$ is nondegenerate on the lowest weight
spaces, and since the left and right radicals are
$\widehat{\mathfrak{g}}$-modules by
invariance, the left radical is contained in $J(\lambda^*)$ and the right
radical
is contained in $J(\lambda)$. On the other hand, suppose $v\in J(\lambda)(n)$.
Since $n>0$, $V^{M(\lambda^*)}(n)$ is spanned by vectors of the form $y\otimes
u$
where $u\in M(\lambda^*)$ and $y\in U(\widehat{\mathfrak{g}}_-)$ has positive
degree. Then
$$ \langle y\otimes u,v\rangle_M =\langle u,\theta(y)v\rangle_M =0$$
because $\theta(y)v$ is in the trivial intersection of $J(\lambda)$ with
$V^{M(\lambda )}(0)$. Thus $v$ is orthogonal to $V^{M(\lambda^*)}(n)$
as well as to every other weight space of $V^{M(\lambda^*)}$ by Proposition
\ref{orthogonal}. This means $J(\lambda)$ is contained in the right radical of
$\langle\cdot,\cdot\rangle_M$. Similarly, $J(\lambda^*)$ is contained in the
left
radical of $\langle\cdot,\cdot\rangle_M$, proving the proposition.
\epfv

Using Proposition \ref{radicals}, we see that $\langle\cdot,\cdot\rangle_M$
induces a well-defined bilinear pairing between any two quotients of
$V^{M(\lambda)}$ and $V^{M(\lambda^*)}$. In particular, the pairing is
nondegenerate
between $L(\lambda)$ and $L(\lambda^*)$. This means that $L(\lambda^*)\cong
L(\lambda)'$. Notice also that the restriction of
$\langle\cdot,\cdot\rangle_M$ to $K(\lambda^*)\times K(\lambda)$ is
nondegenerate.

Let $J(\lambda)^\perp$ denote the set of functionals $u'\in
(V^{M(\lambda)})'$ such that $\langle u',v\rangle =0$ for all $v\in J(\lambda)$. It is easy to see that $J(\lambda)^\perp$ is a $\widehat{\mathfrak{g}}$-module.

\begin{prop}\label{phiprop}
As $\widehat{\mathfrak{g}}$-modules, $J(\lambda)^\perp$ is isomorphic to $L(\lambda^{*})$. Furthermore, there is a $\mathfrak{g}$-module isomorphism $\Phi :
J(\lambda)^\perp\rightarrow K(\lambda^*)$ such that
 \begin{equation*}
  \langle\Phi(u'),v\rangle_M=\langle u',v\rangle
 \end{equation*}
for any $u'\in J(\lambda)^\perp$ and $v\in V^{M(\lambda)}$.
\end{prop}
\pf
 We first define a map $\phi: L(\lambda^*)\rightarrow (V^{M(\lambda)})'$ by
 \begin{equation*}
  \langle\phi(u+J(\lambda^*)),v\rangle =\langle u,v\rangle_M
 \end{equation*}
for $u\in V^{M(\lambda^*)}$ and $v\in V^{M(\lambda)}$. This map is well-defined
and injective because $J(\lambda^*)$ is the left radical of
$\langle\cdot,\cdot\rangle_M$. Furthermore, $\phi$ is a
$\widehat{\mathfrak{g}}$-homomorphism because for $g\in\mathfrak{g}$ and
$n\in\mathbb{Z}$, we have
\begin{align*}
 \langle g(n)\cdot\phi(u+J(\lambda^*)),v\rangle & =
-\langle\phi(u+J(\lambda^*)),g(-n)v\rangle=-\langle u,g(-n)v\rangle_M\nonumber\\
& =  \langle g(n)u,v\rangle_M=\langle\phi(g(n)(u+J(\lambda^*))),v\rangle,
\end{align*}
by the definition of the $\widehat{\mathfrak{g}}$-module action on
$(V^{M(\lambda)})'$. We observe that the image of $\phi$ is contained in
$J(\lambda)^\perp$ because $J(\lambda)$ is the right radical of
$\langle\cdot,\cdot\rangle_M$, and that $\phi$ preserves conformal weights
because it is a $\widehat{\mathfrak{g}}$-homomorphism. Now, the weight spaces of
$J(\lambda)^\perp $ have the same dimension as the weight spaces of
$L(\lambda)$ since $L(\lambda)=V^{M(\lambda)}/J(\lambda)$, and the weight spaces
of $L(\lambda)$ have the same dimension as the weight spaces of $L(\lambda^*)$
since $L(\lambda^*)\cong L(\lambda)'$. Since the weight spaces are finite
dimensional, this means that $\phi$ is surjective onto
$J(\lambda)^\perp$, and we have a $\widehat{\mathfrak{g}}$-isomorphism
$\phi^{-1}: J(\lambda)^\perp\rightarrow L(\lambda^*)$.

We also define a map $\psi: L(\lambda^*)\rightarrow K(\lambda^*)$ by $\psi
(u+J(\lambda^*))=P_K(u)$ for $u\in V^{M(\lambda^*)}$. This is well-defined and
injective because $J(\lambda^*)$ is the kernel of $P_K$ and $\psi$ is a
$\mathfrak{g}$-homomorphism because $P_K$ is a $\mathfrak{g}$-homomorphism.
Also, $\psi$ is surjective because the weight spaces of $L(\lambda^*)$ and
$K(\lambda^*)$ have the same dimensions.

We then define the $\mathfrak{g}$-homomorphism
$\Phi=\psi\circ\phi^{-1}:J(\lambda)^\perp\rightarrow K(\lambda^*)$. If
$u'=\phi(u+J(\lambda^*))\in J(\lambda)^\perp$, then
\begin{equation*}
 \langle \Phi(u'),v\rangle_M =\langle\psi(u+J(\lambda^*)),v\rangle_M =\langle
P_K(u),v\rangle_M =\langle u,v\rangle_M
=\langle\phi(u+J(\lambda^*)),v\rangle=\langle u',v\rangle
\end{equation*}
for any $v\in V^{M(\lambda)}$, as desired.
\epfv
\begin{exam}
 When $\mathfrak{g}=\mathfrak{sl}(2,\mathbb{C})$, we have $\lambda^*=\lambda$
for any
$\lambda=r\alpha/2$ where $r\in\mathbb{N}$. Thus we have constructed a
$\widehat{\mathfrak{g}}$-invariant bilinear form on $V^{M(r)}=V^{M(r\alpha/2)}$
for any
$r\in\mathbb{N}$, and $\Phi$ gives a $\mathfrak{g}$-module isomorphism from
$J(r)^\perp\subseteq (V^{M(r)})'$ to $K(r)$ determined by
\begin{equation}\label{phi}
 \langle\Phi(u'),v\rangle_M =\langle u',v\rangle.
\end{equation}
for any $u'\in J(r)^\perp$, $v\in V^{M(r)}$.
\end{exam}

\setcounter{equation}{0}
\section{The main theorems on intertwining operators among generalized Verma
modules}

\allowdisplaybreaks
In this section we present our main theorems for intertwining operators among generalized Verma modules for $\widehat{\mathfrak{sl}(2,\mathbb{C})}$. We shall be interested in intertwining operators among modules for the
generalized Verma module vertex operator algebra $V^{M(0,\ell)}$ for some fixed $\ell\in\N$. Recalling the notation of \eqref{intwop}, we will prove:
\begin{thm}\label{mainintwopthm}
 Suppose $M(p)$, $M(q)$, and $M(r)$ for $p,q,r\in\N$ are finite-dimensional irreducible $\mathfrak{sl}(2,\C)$-modules with highest weights $p$, $q$, and $r$ respectively and the conditions of Theorem \ref{maintheorem} below hold. Then the linear map
\begin{align*}
\mathcal{V}^{V^{M(r)}}_{V^{M(p)}\,V^{M(q)}} & \longrightarrow\mathrm{Hom}_{\mathfrak{sl}(2,
\mathbb{C})}(M(p)\otimes M(q), M(r))\\
\mathcal{Y} & \longmapsto \mathcal{Y}_0\vert_{M(p)\otimes M(q)}
\end{align*}
is a linear isomorphism.
\end{thm}

In order to prove this theorem, we will need a result on extending an intertwining operator map from lowest weight spaces to entire modules. It is similar to Tsuchiya and Kanie's nuclear democracy theorem \cite{TK} and Li's generalized nuclear democracy theorem (Theorem 4.12 in \cite{Li2}); see also the main theorem in \cite{Li3}. Since the proof requires lengthy but standard calculations using formal calculus, we defer it to the appendix. Note that the following theorem applies to any finite-dimensional complex
simple Lie algebra $\mathfrak{g}$, not merely $\mathfrak{sl}(2,\mathbb{C})$:
\begin{thm}\label{intwopext}
 Suppose $M(\lambda_1)$ and $M(\lambda_2)$ are finite-dimensional irreducible
$\mathfrak{g}$-modules with highest weights $\lambda_1$ and $\lambda_2$, respectively, $W$ is a $V^{M(0)}$-module, and
 \begin{equation*}
  \mathcal{Y}: M(\lambda_1)\otimes M(\lambda_2)\rightarrow
W\lbrace x\rbrace
 \end{equation*}
is a linear map, given by $u\otimes v\mapsto\mathcal{Y}(u,x)v$
for $u\in M(\lambda_1)$ and $v\in M(\lambda_2)$, satisfying
\begin{equation}\label{comm1}
 [g(n),\mathcal{Y}(u,x)]=x^n\mathcal{Y}(g(0)u,x)
\end{equation}
for any $g\in\mathfrak{g}$, $n\geq 0$, and
\begin{equation}\label{l0comm}
 [L(0),\mathcal{Y}(u,x)]=x\dfrac{d}{dx}\mathcal{Y}(u,x)+\mathcal{Y}(L(0)u,x).
\end{equation}
Then $\mathcal{Y}$ has a unique extension to an intertwining operator
\begin{equation*}
 \mathcal{Y}: V^{M(\lambda_1)}\otimes V^{M(\lambda_2)}\rightarrow W\lbrace
x\rbrace.
\end{equation*}
\end{thm}

 Recalling \eqref{intwop}, we see that a map $\mathcal{Y}: M(\lambda_1)\otimes M(\lambda_2)\rightarrow
V^{M(\lambda_3)}\lbrace x\rbrace$ satisfying (\ref{comm1}) and (\ref{l0comm}) is equivalent to a sequence of maps
\begin{equation*}
 \mathcal{Y}_m: M(\lambda_1)\otimes M(\lambda_2)\rightarrow V^{M(\lambda_3)}
\end{equation*}
for each $m\in\Z$ such that $\mathcal{Y}_m(u)v\in V^{M(\lambda_3)}(-m)$ and
\begin{equation}\label{sl2comp1}
 [g(n),\mathcal{Y}_m(u)]v=\mathcal{Y}_{m+n}(g(0)u)v
\end{equation}
for any $u\in M(\lambda_1)$, $v\in M(\lambda_2)$, $m\in\mathbb{Z}$,
$g\in\mathfrak{g}$, and $n\geq 0$. Of course, we take $V^{M(\lambda_3)}(m)=0$ for $m<0$ so that $\mathcal{Y}_m=0$ for $m>0$.

In order to construct maps satisfying \eqref{sl2comp1}, we will take
$\mathfrak{g}=\mathfrak{sl}(2,\mathbb{C})$ once again:
\begin{thm}\label{maintheorem}
Suppose $M(p)$, $M(q)$, and $M(r)$ for $p,q,r\in\N$ are finite-dimensional irreducible $\mathfrak{sl}(2,\C)$-modules with highest weights $p$, $q$, and $r$, respectively, and the maximal proper submodule $J(r)\subseteq V^{M(r)}$ is irreducible
and isomorphic to some $L(r')$. Suppose moreover that $J(r')$ is irreducible and
isomorphic to some $L(r'')$, and that there are no non-zero
$\mathfrak{sl}(2,\mathbb{C})$-homomorphisms from $M(p)\otimes M(q)$ to $M(r')$
or $M(r'')$. Then given
$f\in\mathrm{Hom}_{\mathfrak{sl}(2,\mathbb{C})}(M(p)\otimes M(q),M(r))$, there
are for $m\in\mathbb{Z}$ unique
maps $\mathcal{Y}_m: M(p)\otimes M(q)\rightarrow V^{M(r)}(-m)$ such that $\mathcal{Y}_0=f$ and \eqref{sl2comp1} holds.
\end{thm}

\begin{rema}
 The irreducibility of the maximal proper submodules $J(r)$ and $J(r')$ are natural conditions to consider in light of Theorem \ref{irred}. Moreover, we shall see examples in Section 8 for which there are non-zero homomorphisms $M(p)\otimes M(q)\rightarrow M(r')$ and for which the assertion of Theorem \ref{maintheorem} fails to hold.
\end{rema}

\textit{Proof of Theorem 6.3.}\hspace{1ex} We take $\mathfrak{g}=\mathfrak{sl}(2,\C)$-modules $M(p)$ and $M(q)$ and a
generalized Verma module $V^{M(r)}$. We assume that $J(r)$ is irreducible and
isomorphic to some $L(r')$ and that $J(r')$ is irreducible and isomorphic to
some $L(r'')$. Since $J(r)\cong L(r')$, we may equip $J(r)$ with the
nondegenerate invariant bilinear form $\langle\cdot,\cdot\rangle_J$ induced by
the invariant form on $V^{M(r')}$ constructed in the previous section. Suppose
the lowest weight space of $J(r)$ is contained in $V^{M(r)}(M)$ for some $M>0$.
Note that in the identification $J(r)\cong L(r')$, the weight space $J(r)(n)$ for any $n\geq M$ is
then identified with $L(r')(n-M)$ (recall Remark \ref{gradingnotation}).

In order to prove the theorem, take a $\mathfrak{g}$-homomorphism
$f: M(p)\otimes M(q)\rightarrow M(r)$. We need to construct unique maps
\begin{equation*}
 \mathcal{Y}_m: M(p)\otimes M(q)\rightarrow V^{M(r)}(-m)
\end{equation*}
for all $m\in\mathbb{Z}$ such that for any $u\in M(p)$, $v\in M(q)$,
\begin{equation}\label{base}
 \mathcal{Y}_0(u)v=f(u\otimes v),
\end{equation}
and
\begin{equation}\label{sl2comp3}
  [g(n),\mathcal{Y}_m(u)]v=\mathcal{Y}_{m+n}(g(0)u)v
\end{equation}
for $g\in\mathfrak{g}$ and $n\geq 0$.

 For the
uniqueness assertion, suppose we have maps $\mathcal{Y}_m$ satisfying
(\ref{base}) and (\ref{sl2comp3}). We write
$\mathcal{Y}_m=\mathcal{Y}^K_m+\mathcal{Y}^J_m$ where $\mathcal{Y}^K_m=P_K\circ
\mathcal{Y}_m$ and $\mathcal{Y}^J_m=P_J\circ \mathcal{Y}_m$.  Then it is clear
that $\mathcal{Y}^K_m=0$ for $m>0$ and $\mathcal{Y}^J_m=0$ for $m>-M$, and that
$\mathcal{Y}^K_0=f$. We also claim that $\mathcal{Y}^J_{-M}=0$. Indeed, by the $n=0$ case of (\ref{sl2comp3}), each
$\mathcal{Y}_m$ is a $\mathfrak{g}$-module homomorphism from
$M(p)\otimes M(q)$ to
$V^{M(r)}(-m)$, and because $P_K$ and $P_J$ are $\mathfrak{g}$-module
homomorphisms,
so is each $\mathcal{Y}^K_m$ and $\mathcal{Y}^J_m$. In particular,
$\mathcal{Y}^J_{-M}$ is a $\mathfrak{g}$-module homomorphism from $M(p)\otimes
M(q)$ to $M(r')$, the lowest weight space of $J(r)$. However, by the assumptions
of the theorem, the only such homomorphism is $0$.

 Now,  for any $m\in\mathbb{Z}$, $n>0$, $g\in\mathfrak{g}$, $u\in M(p)$,
and $v\in M(q)$, (\ref{sl2comp3}) implies that
 \begin{align}\label{split}
\mathcal{Y}^K_{-m+n}(g(0)u)v+ & \mathcal{Y}^J_{-m+n}(g(0)u)v =
\mathcal{Y}_{-m+n}(g(0)u)v =g(n)\mathcal{Y}_{-m}(u)v \nonumber\\
& =g(n)\mathcal{Y}^K_{-m}(u)v+g(n)\mathcal{Y}^J_{-m}(u)v\nonumber\\
&
=P_K(g(n)\mathcal{Y}^K_{-m}(u)v)+P_J(g(n)\mathcal{Y}^K_{-m}(u)v)+g(n)\mathcal{Y}
^J_{-m}
(u)v.
 \end{align}
Thus if we suppose $m>0$, $0<n\leq m$ and $w\in V^{M(r)}(m-n)$, the invariance
of the form $\langle\cdot,\cdot\rangle_M$ on $V^{M(r)}$, the fact that $J(r)$ is
the radical of the form, and the projection of (\ref{split}) to $K(r)$ imply
that
\begin{align*}
 \langle \mathcal{Y}^K_{-m}(u)v, g(-n)w\rangle_M & =  -\langle
g(n)\mathcal{Y}^K_{-m}(u)v,w\rangle_M\nn
& = -\langle P_K(g(n)\mathcal{Y}^K_{-m}(u)v),w\rangle_M\nn
& = -\langle\mathcal{Y}^K_{-m+n}(g(0)u)v,w\rangle_M.
\end{align*}
Since $\langle\cdot,\cdot\rangle_M$ is nondegenerate on $K(r)$ and since
$K(r)(m)$ is
spanned by linear combinations of certain vectors of the form $g(-n)w$ where
$0<n\leq m$ and $w\in V^{M(r)}(m-n)$,  we
see that for $m>0$, $\mathcal{Y}^K_{-m}$ is uniquely determined by
$\mathcal{Y}^K_{-m+n}$ for $0<n\leq m$. Thus $\mathcal{Y}^K_{-m}$ for $m>0$ is
uniquely determined by $\mathcal{Y}^K_0=f$.

Similarly, if we suppose $m>M$, $0<n\leq m-M$, and $w\in J(r)(m-n)$, we must
have
\begin{align*}
  \langle \mathcal{Y}^J_{-m}(u)v, g(-n)w\rangle_J & =  -\langle
g(n)\mathcal{Y}^J_{-m}(u)v,w\rangle_J\nonumber\\
  & =  -\langle \mathcal{Y}^J_{-m+n}(g(0)u)v,w\rangle_J+\langle
P_J(g(n)\mathcal{Y}^K_{-m}(u)v),w\rangle_J.
\end{align*}
Since $\langle\cdot,\cdot\rangle_J$ is nondegenerate on $J(r)$ and since
vectors of the form $g(-n)w$ for $0<n\leq m-M$ and $w\in J(r)(m-n)$ span
$J(r)(m)$,
we see as before that $\mathcal{Y}^J_{-m}$ for $m>M$ is determined by
$\mathcal{Y}^J_{-M}=0$ and by $\mathcal{Y}^K_{-m}$. This proves the uniqueness
assertion of the theorem, and it remains to show that the above recursive
formulas
for $\mathcal{Y}^K_{-m}$ and $\mathcal{Y}^J_{-m}$ do in fact determine
well-defined linear maps with the required properties.

 Since $M(r)\cong M(r)^*$, the lowest weight space of $(V^{M(r)})'$ is
isomorphic to $M(r)$. Thus we can view $f$ as a homomorphism from $M(p)\otimes
M(q)$ into $(V^{M(r)})'(0)$. Now we define
\begin{equation*}
 \bar{\mathcal{Y}}^K_{-m}: M(p)\otimes M(q)\rightarrow (V^{M(r)})'(m)
\end{equation*}
for any $m\geq 0$ recursively using the formulas:
\begin{equation}\label{basecase}
 \langle \bar{\mathcal{Y}}^K_0(u)v,w\rangle=\langle f(u\otimes v),w\rangle
\end{equation}
for $u\in M(p)$, $v\in M(q)$, and $w\in M(r)$, and for $m>0$,
\begin{equation}\label{recursion}
 \langle \bar{\mathcal{Y}}^K_{-m}(u)v,g(-n)w\rangle =-\langle
\bar{\mathcal{Y}}^K_{-m+n}(g(0)u)v,w\rangle,
\end{equation}
where $g\in\mathfrak{g}$, $0<n\leq m$, and $w\in V^{M(r)}(m-n)$. Note that
$V^{M(r)}(m)$ is spanned by vectors of the form $g(-n)w$, so \eqref{recursion}
does define a vector in $(V^{M(r)})'(m)$.

\begin{lemma}
 For any $m\geq 0$, we have
 \begin{equation*}
  [g(0),\bar{\mathcal{Y}}^K_{-m}(u)]v=\bar{\mathcal{Y}}^K_{-m}(g(0)u)v
 \end{equation*}
for $g\in\mathfrak{g}$, $u\in M(p)$, and $v\in M(q)$.
\end{lemma}
\pf
We prove by induction on $m$. The conclusion is true for $m=0$ because
$\bar{\mathcal{Y}}^K_{0}=f$ is a $\mathfrak{g}$-module homomorphism. For $m>0$,
we have
\begin{align*}
 \langle  [g(0), & \bar{\mathcal{Y}}^K_{-m}(u)]v, h(-n)w\rangle  =  -\langle
\bar{\mathcal{Y}}^K_{-m}(u)v, g(0)h(-n)w\rangle -\langle
\bar{\mathcal{Y}}^K_{-m}(u)g(0)v, h(-n)w\rangle\nonumber\\
 & =  -\langle \bar{\mathcal{Y}}^K_{-m}(u)v, h(-n)g(0)w\rangle-\langle
\bar{\mathcal{Y}}^K_{-m}(u)v,[g,h](-n)w\rangle+\langle
\bar{\mathcal{Y}}^K_{-m+n}(h(0)u)g(0)v, w\rangle\nonumber\\
 & =  -\langle  g(0)\bar{\mathcal{Y}}^K_{-m+n}(h(0)u)v , w\rangle +\langle
\bar{\mathcal{Y}}^K_{-m+n}([g,h](0)u)v, w\rangle+\langle
\bar{\mathcal{Y}}^K_{-m+n}(h(0)u)g(0)v,w\rangle\nonumber\\
 & =  -\langle [g(0), \bar{\mathcal{Y}}^K_{-m+n}(h(0)u)]v, w\rangle +\langle
\bar{\mathcal{Y}}^K_{-m+n}([g(0),h(0)]u)v ,w\rangle\nonumber\\
 & =  -\langle  \bar{\mathcal{Y}}^K_{-m+n}(g(0)h(0)u)v, w\rangle +\langle
\bar{\mathcal{Y}}^K_{-m+n}(g(0)h(0)u)v ,w\rangle  -\langle
\bar{\mathcal{Y}}^K_{-m+n}(h(0)g(0)u)v, w\rangle \nonumber\\
 & =  \langle  \bar{\mathcal{Y}}^K_{-m}(g(0)u)v , h(-n)w\rangle .
\end{align*}
Since $V^{M(r)}(m)$ is spanned by vectors of the form $h(-n)w$, for
$h\in\mathfrak{g}$, $0<n\leq m$ and $w\in V^{M(r)}(m-n)$, the conclusion
follows.
\epfv
\begin{lemma}\label{radlem}
 For any $m\geq 0$, $u\in M(p)$, and $v\in M(q)$, $
\bar{\mathcal{Y}}^K_{-m}(u)v\in J(r)^\perp$.
\end{lemma}
\pf
By the previous lemma, $\bar{\mathcal{Y}}^K_{-m}$ is a $\mathfrak{g}$-module
homomorphism from $M(p)\otimes M(q)$ to $(V^{M(r)})'(m)$. Also, as a
$\mathfrak{g}$-module, the lowest weight space of $J(r)\cong L(r')$ is
$M(r')$. Thus $\bar{\mathcal{Y}}^K_{-m}$ defines a $\mathfrak{g}$-module
homomorphism $f_m$ from $M(p)\otimes M(q)$ to $M(r')^*\cong M(r')$ via
\begin{equation*}
\langle f_m(u\otimes v),w\rangle = \langle\bar{\mathcal{Y}}^K_{-m}(u)v,w\rangle
\end{equation*}
for $u\in M(p)$, $v\in M(q)$, and $w\in M(r')$. But the assumptions of Theorem
\ref{maintheorem} imply that there are
no non-zero $\mathfrak{g}$-module homomorphisms from $M(p)\otimes M(q)$ to
$M(r')$. Thus
\begin{equation*}
 \langle\bar{\mathcal{Y}}^K_{-m}(u)v,w\rangle =0
\end{equation*}
for any $w\in M(r')=J(r')(M)$.

Now, if $w\in J(r)$ is a homogeneous vector such that
$\langle\bar{\mathcal{Y}}^K_{-m}(u)v,w\rangle =0$ for any $u\in M(p)$, $v\in
M(q)$, and $m\geq 0$, then for any $g\in\mathfrak{g}$ and $n>0$ such that
$g(-n)w\in J(r)(m)$, we have
\begin{equation*}
\langle\bar{\mathcal{Y}}^K_{-m}(u)v,g(-n)w\rangle
=-\langle\bar{\mathcal{Y}}^K_{-m+n}(g(0)u)v,w\rangle=0.
\end{equation*}
If $n>0$ and $g(-n)w\notin J(r)(m)$, then
$\langle\bar{\mathcal{Y}}^K_{-m}(u)v,g(-n)w\rangle =0$ automatically. Since
$J(r)\cong L(r')$, $J(r)$ is generated as $\widehat{\mathfrak{g}}_-$-module by
its lowest weight space $M(r')$. Thus
$\langle\bar{\mathcal{Y}}^K_{-m}(u)v,w\rangle =0$ for all $w\in J(r)$, and
$\bar{\mathcal{Y}}^K_{-m}(u)v\in J(r)^\perp$.
\epfv

Now, by the previous lemma, we can use the $\mathfrak{g}$-module isomorphism
$\Phi: J(r)^\perp\rightarrow K(r)$ given by Proposition \ref{phiprop} to
identify
$\bar{\mathcal{Y}}^K_{-m}(u)v$ with
$\mathcal{Y}^K_{-m}(u)v=\Phi(\bar{\mathcal{Y}}^K_{-m}(u)v)\in K(r)$. By
(\ref{phi}), (\ref{basecase}), and (\ref{recursion}), we have that
$\mathcal{Y}^K_{0}(u)v=f(u\otimes v)$ and for $m>0$,
\begin{equation*}
 \langle\mathcal{Y}^K_{-m}(u)v , g(-n)w\rangle_M =-\langle
\mathcal{Y}^K_{-m+n}(g(0)u)v,w\rangle_M
\end{equation*}
for $u\in M(p)$, $v\in M(q)$ and $g(-n)w\in V^{M(r)}(m)$. Moreover, since $\Phi$
is a
$\mathfrak{g}$-module homomorphism, we also have for $m\geq 0$
 \begin{equation}\label{khom}
  [g(0),\mathcal{Y}^K_{-m}(u)]v=\mathcal{Y}^K_{-m}(g(0)u)v
 \end{equation}
for $g\in\mathfrak{g}$, $u\in M(p)$, and $v\in M(q)$.

For $m<0$, we set $\mathcal{Y}^K_{-m}=0$. Then we have:
\begin{lemma}\label{khom2}
 If $n>0$, $P_K(g(n)\mathcal{Y}^K_{-m}(u)v)=\mathcal{Y}^K_{-m+n}(g(0)u)v$ for
$u\in M(p)$, $v\in M(q)$, $g\in\mathfrak{g}$, and any $m\in\mathbb{Z}$.
\end{lemma}
\pf
If $n>m$, both sides are zero. Otherwise, if $0<n\leq m$ (so that $m>0$), for
any
$w\in V^{M(r)}(m-n)$ we have
\begin{align*}
 \langle P_K(g(n)\mathcal{Y}^K_{-m}(u)v),w\rangle_M & =  \langle
g(n)\mathcal{Y}^K_{-m}(u)v,w\rangle_M\nonumber\\
 & =  -\langle \mathcal{Y}^K_{-m}(u)v,g(-n)w\rangle_M\nonumber\\
 & =  \langle \mathcal{Y}^K_{-m+n}(g(0)u)v,w\rangle_M.
\end{align*}
Since $\langle\cdot,\cdot\rangle_M$ is nondegenerate on $K(r)$, and since both
$P_K(g(n)\mathcal{Y}^K_{-m}(u)v)$ and
$\mathcal{Y}^K_{-m+n}(g(0)u)v$ are in $K(r)(m-n)$, they are equal.
\epfv

Next, we want to construct maps $\bar{\mathcal{Y}}^J_{-m}: M(p)\otimes
M(q)\rightarrow (V^{M(r')})'(m-M)$ for $m\geq M$. We define these maps
recursively:
$\bar{\mathcal{Y}}^J_{-M}=0$ and for $m>M$ we define
\begin{equation*}
 \langle\bar{\mathcal{Y}}^J_{-m}(u)v, g(-n)w\rangle =\langle
P_J(g(n)\mathcal{Y}^K_{-m}(u)v),w+J(r')\rangle_J -\langle
\bar{\mathcal{Y}}^J_{-m+n}(g(0)u)v,w\rangle
\end{equation*}
for $u\in M(p)$, $v\in M(q)$, $g\in\mathfrak{g}$, $0<n\leq m-M$, and $w\in
V^{M(r')}(m-M-n)$. Notice that in the first term on the right, we are
identifying
the quotient $V^{M(r')}/J(r')$ with $J(r)\subseteq V^{M(r)}$.
\begin{rema}
 Since we will need $\bar{\mathcal{Y}}^J_{-M}$ to be a
$\mathfrak{g}$-homomorphism from $M(p)\otimes M(q)$ into $(V^{M(r')})'(0)\cong
M(r')^*\cong M(r')$, and since by assumption there are no non-zero such
homomorphisms, we are required to set $\bar{\mathcal{Y}}^J_{-M}=0$.
\end{rema}
\begin{lemma}\label{jhom}
 For any $m\geq M$, we have
  \begin{equation*}
  [g(0),\bar{\mathcal{Y}}^J_{-m}(u)]v=\bar{\mathcal{Y}}^J_{-m}(g(0)u)v
 \end{equation*}
for $g\in\mathfrak{g}$, $u\in M(p)$, and $v\in M(q)$.
\end{lemma}
\pf
We prove by induction on $m$. For $m=M$, the result is true because both sides
of the equation are zero. Now for $m>M$, using the induction hypothesis, the
fact that $P_J$ is a $\mathfrak{g}$-module homomorphism, and (\ref{khom}), we
have for $h\in\mathfrak{g}$, $0<n\leq m-M$, and $w\in V^{M(r')}(m-M-n)$:
\begin{align*}
 \langle [g(0),\bar{\mathcal{Y}}^J_{-m} & (u)]v,\, h(- n)w\rangle  =  -\langle
\bar{\mathcal{Y}}^J_{-m}(u)v, g(0)h(-n)w\rangle -\langle
\bar{\mathcal{Y}}^J_{-m}(u)g(0)v, h(-n)w\rangle\nonumber\\
 & =  -\langle \bar{\mathcal{Y}}^J_{-m}(u)v, h(-n)g(0)w\rangle -\langle
\bar{\mathcal{Y}}^J_{-m}(u)v, [g,h](-n)w\rangle\nonumber\\
 &  \hspace{1em} +\langle
\bar{\mathcal{Y}}^J_{-m+n}(h(0)u)g(0)v,w\rangle-\langle
P_J(h(n)\mathcal{Y}^K_{-m}(u)g(0)v),w+J(r')\rangle_J\nonumber\\
 & =  -\langle g(0)\bar{\mathcal{Y}}^J_{-m+n}(h(0)u)v,w\rangle+\langle
g(0)P_J(h(n)\mathcal{Y}^K_{-m}(u)v), w+J(r')\rangle_J\nonumber\\
 & \hspace{1em} +\langle \bar{\mathcal{Y}}^J_{-m+n}([g,h](0)u)v,w\rangle-\langle
P_J([g,h](n)\mathcal{Y}^K_{-m}(u)v),w+J(r')\rangle_J\nonumber\\
 & \hspace{1em} +\langle \bar{\mathcal{Y}}^J_{-m+n}(h(0)u)g(0)v,w\rangle-\langle
P_J(h(n)\mathcal{Y}^K_{-m}(u)g(0)v),w+J(r')\rangle_J\nonumber\\
 & =  -\langle [g(0),\bar{\mathcal{Y}}^J_{-m+n}(h(0)u)]v,w\rangle +\langle
\bar{\mathcal{Y}}^J_{-m+n}([g(0),h(0)]u)v,w\rangle\nonumber\\
 & \hspace{1em} +\langle
P_J(g(0)h(n)\mathcal{Y}^K_{-m}(u)v-[g(0),h(n)]\mathcal{Y}^K_{-m}(u)v),
w+J(r')\rangle_J\nonumber\\
 & \hspace{1em} -\langle
P_J(h(n)\mathcal{Y}^K_{-m}(u)g(0)v),w+J(r')\rangle_J\nonumber\\
 & =  -\langle \bar{\mathcal{Y}}^J_{-m+n}(g(0)h(0)u)v,w\rangle +\langle
\bar{\mathcal{Y}}^J_{-m+n}(g(0)h(0)u)v,w\rangle\nonumber\\
 & \hspace{1em} -\langle \bar{\mathcal{Y}}^J_{-m+n}(h(0)g(0)u)v,w\rangle
+\langle
P_J(h(n)[g(0),\mathcal{Y}^K_{-m}(u)]v),w+J(r')\rangle_J\nonumber\\
 & =  -\langle \bar{\mathcal{Y}}^J_{-m+n}(h(0)g(0)u)v,w\rangle +\langle
P_J(h(n)\mathcal{Y}^K_{-m}(g(0)u)v),w+J(r')\rangle_J\nonumber\\
 & =  \langle \bar{\mathcal{Y}}^J_{-m}(g(0)u)v,h(-n)w\rangle.
\end{align*}
Since $V^{M(r')}(m-M)$ is spanned by the vectors $h(-n)w$ where
$h\in\mathfrak{g}$, $0<n\leq m-M$ and $w\in V^{M(r')}(m-M-n)$,  the result
follows.
\epfv
\begin{lemma}
 For any $m\geq M$, $u\in M(p)$, and $v\in M(q)$, $
\bar{\mathcal{Y}}^J_{-m}(u)v\in J(r')^\perp$.
\end{lemma}
\pf
The previous lemma implies that $ \bar{\mathcal{Y}}^J_{-m}$ determines a
$\mathfrak{g}$-module homomorphism from $M(p)\otimes M(q)$ to $M(r'')^*\cong
M(r'')$, where $M(r'')$ is the lowest weight space of
$J(r')\cong L(r'')$; this homomorphism is necessarily zero, just as in the proof
of Lemma
\ref{radlem}. Thus
\begin{equation*}
 \langle\bar{\mathcal{Y}}^J_{-m}(u)v,w\rangle =0
\end{equation*}
for any $w\in M(r'')\subseteq J(r')$.

If $w\in J(r')$ is a homogeneous vector such that
$\langle\bar{\mathcal{Y}}^J_{-m}(u)v,w\rangle =0$ for any $u\in M(p)$, $v\in
M(q)$, and $m\geq M$, then for any $g\in\mathfrak{g}$ and $n>0$ such that
$g(-n)w\in J(r')(m-M)$,
\begin{equation*}
\langle\bar{\mathcal{Y}}^J_{-m}(u)v,g(-n)w\rangle
=-\langle\bar{\mathcal{Y}}^J_{-m+n}(g(0)u)v,w\rangle+\langle
P_J(g(n)\mathcal{Y}^K_{-m}(u)v),w+J(r')\rangle_J=0
\end{equation*}
because $w+J(r')=0$. If $g(-n)w\notin J(r')(m-M)$, then
$\langle\bar{\mathcal{Y}}^J_{-m}(u)v,g(-n)w\rangle =0$ automatically.  As in
the
proof of Lemma \ref{radlem}, this proves the lemma.
\epfv

By this last lemma, $\bar{\mathcal{Y}}^J_{-m}(u)v$ for $m\geq M$ defines an
element in $(V^{M(r')}/J(r'))'\cong L(r')'\cong J(r)'$. Since the nondegenerate
form
$\langle\cdot,\cdot\rangle_J$ identifies $J(r)'$ with $J(r)$,
$\bar{\mathcal{Y}}^J_{-m}(u)v$ induces a unique element
$\mathcal{Y}^J_{-m}(u)v\in J(r)(m)\cong
L(r')(m-M)$ such that $\mathcal{Y}^J_{-M}(u)v=0$ and for $m>M$
\begin{equation*}
  \langle\mathcal{Y}^J_{-m}(u)v, g(-n)w\rangle_J =\langle
P_J(g(n)\mathcal{Y}^K_{-m}(u)v),w\rangle_J -\langle
\mathcal{Y}^J_{-m+n}(g(0)u)v,w\rangle_J
\end{equation*}
for any $g\in\mathfrak{g}$, $0<n\leq m-M$, and $w\in J(r)(m-n)$; also
\begin{equation}\label{jcomp}
  [g(0),\mathcal{Y}^J_{-m}(u)]v=\mathcal{Y}^J_{-m}(g(0)u)v.
\end{equation}
We define $\mathcal{Y}^J_{-m}(u)v=0$ for $m<M$, and so we have:
\begin{lemma}\label{jcomp2}
 For any $m\in\mathbb{Z}$ and  $n>0$,
 \begin{equation*}
P_J(g(n)\mathcal{Y}^K_{-m}(u)v)+g(n)\mathcal{Y}^J_{-m}(u)v=\mathcal{Y}^J_{-m+n}
(g(0)u)v,
 \end{equation*}
where $g\in\mathfrak{g}$, $u\in M(p)$, and $v\in M(q)$.
\end{lemma}
\pf
If $n>m-M$, both sides are zero. Otherwise, if $0<n\leq m-M$ (so that $m>M$),
then
for any $w\in J(r)(m-n)$,
\begin{align*}
 \langle g(n)\mathcal{Y}^J_{-m}(u)v,w\rangle_J & = -\langle
\mathcal{Y}^J_{-m}(u)v, g(-n)w\rangle_J\nonumber\\
 & = \langle \mathcal{Y}^J_{-m+n}(g(0)u)v,w\rangle_J-\langle
P_J(g(n)\mathcal{Y}^K_{-m}(u)v),w\rangle_J.
\end{align*}
Since $\langle\cdot,\cdot\rangle_J$ is nondegenerate on $J(r)(m-n)$, the result
follows.
\epfv

Finally, we can define
$\mathcal{Y}_m(u)v=\mathcal{Y}^K_{m}(u)v+\mathcal{Y}^J_{m}(u)v$ for any
$m\in\mathbb{Z}$. We have $\mathcal{Y}_0(u)v=f(u\otimes v)$, and
\begin{equation}
   [g(0),\mathcal{Y}_m(u)]v=\mathcal{Y}_{m}(g(0)u)v
\end{equation}
by (\ref{khom}) and (\ref{jcomp}) (or by the fact that both sides are zero if
$m>0$).
Also, by Lemmas \ref{khom2} and \ref{jcomp2}, for $n>0$,
\begin{align*}
 [g(n),\mathcal{Y}_m(u)]v & = g(n)\mathcal{Y}_m(u)v\nonumber\\
 & = g(n)\mathcal{Y}^K_m(u)v+g(n)\mathcal{Y}^J_m(u)v\nonumber\\
 & =
P_K(g(n)\mathcal{Y}^K_m(u)v)+P_J(g(n)\mathcal{Y}^K_m(u)v)+g(n)\mathcal{Y}
^J_m(u)v\nonumber\\
 & = \mathcal{Y}^K_{m+n}(g(0)u)v+\mathcal{Y}^J_{m+n}(g(0)u)v\nonumber\\
 & = \mathcal{Y}_{m+n}(g(0)u)v.
\end{align*}
Thus we have proven the existence of the desired maps $\mathcal{Y}_m$,
completing the proof of the theorem.
\epfv

Now using Theorems \ref{intwopext} and \ref{maintheorem}, we can complete the proof of Theorem \ref{mainintwopthm}:

\pf If $\mathcal{Y}\in \mathcal{V}^{V^{M(r)}}_{V^{M(p)}\,V^{M(q)}}$, then
(\ref{sl2comp1}) implies that $\mathcal{Y}_0\vert_{M(p)\otimes M(q)}$ gives an
$\mathfrak{sl}(2,\mathbb{C})$-module homomorphism $f_\mathcal{Y}$ into $M(r)$.
Conversely, Theorems \ref{intwopext} and \ref{maintheorem} give us a unique
intertwining operator $\mathcal{Y}_f$ for each
$f\in\mathrm{Hom}_{\mathfrak{sl}(2,\mathbb{C})}(M(p)\otimes M(q),M(r))$ such that $(\mathcal{Y}_f)_0(u)v=f(u\otimes v)$ for $u\otimes
v\in M(p)\otimes M(q)$. So $f_{\mathcal{Y}_f}=f$, and the uniqueness
assertions in Theorems \ref{intwopext} and \ref{maintheorem} guarantee that
$\mathcal{Y}_{f_\mathcal{Y}}=\mathcal{Y}$. Thus we have a linear isomorphism
between $\mathcal{V}^{V^{M(r)}}_{V^{M(p)}\,V^{M(q)}}$ and
$\mathrm{Hom}_{\mathfrak{sl}(2,\mathbb{C})}(M(p)\otimes M(q),M(r))$.
\epfv

\setcounter{equation}{0}
\section{
Intertwining operators of type $\binom{L(r)}{V^{M(p)}\,L(q)}$}

In this section we continue to take $\mathfrak{g}=\mathfrak{sl}(2,\mathbb{C})$
and continue to fix a level $\ell\in\N$. We prove that under certain circumstances,
an intertwining operator of type $\binom{V^{M(r)}}{V^{M(p)}\,V^{M(q)}}$ descends
to an intertwining operator of type $\binom{L(r)}{V^{M(p)}\,L(q)}$.
\begin{thm}\label{maintheorem2}
 Assume the conditions of Theorem \ref{maintheorem} hold and that in addition
$J(q)$ is generated by its lowest weight space $M(q')$ and that there are no
non-zero $\mathfrak{g}$-homomorphisms from $M(p)\otimes M(r)$ to $M(q')$. Then there is a natural linear isomorphism
$\mathcal{V}^{L(r)}_{V^{M(p)}\,L(q)}\cong
\mathcal{V}^{V^{M(r)}}_{V^{M(p)}\,V^{M(q)}}$.
\end{thm}
\pf
We want to show that if $\mathcal{Y}$ is an intertwining operator of type
$\binom{ V^{M(r)}}{V^{M(p)}\,V^{M(q)}}$, then we obtain a well-defined
intertwining operator $\mathcal{Y}^L$ of type $\binom{L(r)}
{V^{M(p)}\,L(q)}$ defined by
\begin{equation}\label{YL}
 \mathcal{Y}^L(u,x)(v+J(q))=\mathcal{Y}(u,x)v+J(r).
\end{equation}
Then the linear map $\mathcal{Y}\mapsto\mathcal{Y}^L$ will be our desired
isormorphism. It will be injective because both $\mathcal{Y}$ and
$\mathcal{Y}^L$ define the same $\mathfrak{g}$-homomorphism
\begin{equation*}
 \mathcal{Y}_0=\mathcal{Y}^L_0: M(p)\otimes M(q)\rightarrow M(r),
\end{equation*}
and by Theorem \ref{mainintwopthm}, $\mathcal{Y}$ is uniquely determined by this
homomorphism. Moreover, the map will be surjective because an argument similar
to but simpler than the uniqueness argument in the proof of Theorem
\ref{maintheorem} shows that any intertwining operator $\mathcal{Y}^L$ of type
$\binom{L(r)}{V^{M(p)}\,L(q)}$ is also determined by such a homomorphism
$M(p)\otimes M(q)\rightarrow M(r)$, and so the fusion rule
$\mathcal{N}^{L(r)}_{V^{M(p)}\,L(q)}$ is  no more than
$\mathcal{N}^{V^{M(r)}}_{V^{M(p)}\,V^{M(q)}}$.

Because $\mathcal{Y}$ is an intertwining operator, $\mathcal{Y}^L$ will satisfy
lower truncation, the Jacobi identity, and the $L(-1)$-derivative property,
provided that it is well defined. To prove well-definedness, we need to show
that if $u\in V^{M(p)}$ and $v\in J(q)$, then $\mathcal{Y}(u,x)v\in J(r)\lbrace
x\rbrace$.
\begin{lemma}
 The subspace $W$ of $V^{M(r)}$ spanned by the coefficients of products of the
form $\mathcal{Y}(u,x)v$ where $u\in V^{M(p)}$ and $v\in J(q)$ is a
$\widehat{\mathfrak{g}}$-submodule of $V^{M(r)}$.
 \end{lemma}
 \pf
  If $g\in\mathfrak{g}$, $n\in\mathbb{Z}$, $u\in V^{M(p)}$, and $v\in J(q)$,
then from the commutator formula,
  \begin{equation*}
   g(n)\mathcal{Y}(u,x)v=\mathcal{Y}(u,x)g(n)v+\sum_{i\geq 0} \binom{n}{i}
x^{n-i}\mathcal{Y}(g(i)u,x)v.
  \end{equation*}
Since all coefficients on the right side are in $W$, the result follows.
 \epfv
 \begin{lemma}
  The intersection of $W$ with the lowest weight space of $V^{M(r)}$ is trivial.
 \end{lemma}
 \pf
  If we use $\langle\cdot,\cdot\rangle_M$ to denote a non-zero
$\widehat{\mathfrak{g}}$-invariant bilinear form on $V^{M(r)}$ as in Proposition
\ref{invariantform}, then it is enough to show that
  \begin{equation}\label{nointersection}
   \langle\mathcal{Y}(u,x)v,w\rangle_M =0
  \end{equation}
for any $u\in V^{M(p)}$, $v\in J(q)$, and $w\in M(r)$, since
$\langle\cdot,\cdot\rangle_M$ gives a nondegenerate bilinear form on $M(r)$. We
first assume $u\in M(p)$ and $v\in M(q')$ (the lowest weight space of $J(q)$).
Then for each $n\in\mathbb{C}$, we have a map $F_n: M(p)\otimes M(r)\rightarrow
M(q')^*\cong M(q')$ given by
\begin{equation*}
 \langle F_n(u\otimes w),v\rangle =\langle\mathcal{Y}_n(u)v,w\rangle_M
\end{equation*}
for $u\in M(p)$, $v\in M(q')$, and $w\in M(r)$. These maps are
$\mathfrak{g}$-homomorphisms because if $g\in\mathfrak{g}$,
\begin{align*}
 \langle g\cdot F_n(u\otimes w),v\rangle & =  -\langle F_n(u\otimes w),
g(0)v\rangle =-\langle\mathcal{Y}_n(u)g(0)v,w\rangle_M\nonumber\\
 & =  -\langle g(0)\mathcal{Y}_n(u)v,w\rangle_M
+\langle\mathcal{Y}_n(g(0)u)v,w\rangle_M\nonumber\\
 & =
\langle\mathcal{Y}_n(u)v,g(0)w\rangle_M+\langle\mathcal{Y}_n(g(0)u)v,
w\rangle_M\nonumber\\
 & =  \langle F_n(g(0)u\otimes w+u\otimes g(0)w),v\rangle.
\end{align*}
But by assumption there are no non-zero homomorphisms from $M(p)\otimes M(r)$ to
$M(q')$, so therefore
\begin{equation*}
 \langle\mathcal{Y}(u,x)v,w\rangle_M =0
\end{equation*}
when $u\in M(p)$ and $v\in M(q')$. Now we suppose that (\ref{nointersection})
holds when $u\in M(p)$ and $w\in M(r)$ for some $v\in J(q)$ and show that it
also holds for $g(-n)v$ where $g\in\mathfrak{g}$ and $n>0$; this will show that
(\ref{nointersection}) holds for all $v\in J(q)$ since $J(q)$ is generated as
$\widehat{\mathfrak{g}}_-$-module by its lowest weight space. Thus:
\begin{equation*}
 \langle\mathcal{Y}(u,x)g(-n)v,w\rangle_M = \langle
g(-n)\mathcal{Y}(u,x)v,w\rangle_M-\langle \mathcal{Y}(g(0)u,x)v,w\rangle_M =0
\end{equation*}
since $g(-n)\mathcal{Y}(u,x)v$ has no component in the lowest weight space since
$n>0$. Finally, we must prove that if (\ref{nointersection}) holds for any $v\in
J(q)$, $w\in M(r)$, and for some $u\in V^{M(p)}$, then it also holds for
$g(-n)u$ when $g\in\mathfrak{g}$ and $n>0$:
\begin{align*}
 \langle\mathcal{Y}(g(-n)u,x)v,w\rangle_M & =
\mathrm{Res}_{x_1}(x_1-x)^{-n}\langle
g(x_1)\mathcal{Y}(u,x)v,w\rangle_M\nonumber\\
&\hspace{1em}-\mathrm{Res}_{x_1} (-x+x_1)^{-n}\langle
\mathcal{Y}(u,x)g(x_1)v,w\rangle_M\nonumber\\
 & =  \sum_{i\geq 0} \binom{-n}{i} (-x)^i\langle
g(-n-i)\mathcal{Y}(u,x)v,w\rangle_M +0 = 0
\end{align*}
by the assumptions on $u$ and since $n+i>0$.
\epfv

Thus we have shown that $W$ is a submodule of $V^{M(r)}$ such that $W\cap
M(r)=0$. This means that $W$ is a proper submodule contained in $J(r)$, and so
$\mathcal{Y}(u,x)v\subseteq J(r)\lbrace x\rbrace $ for any $u\in V^{M(p)}$,
$v\in J(q)$, completing the proof of the theorem.
\epfv

\begin{rema}
 Applying skew-symmetry (see for example \cite{FHL}) to Theorem
\ref{maintheorem2} immediately yields a corresponding theorem on intertwining operators of
type $\binom{L(r)}{L(p)\,V^{M(q)}}$.
\end{rema}

\setcounter{equation}{0}
\section{Examples and counterexamples}

In this section we fix a level $\ell\in\N$ for $\widehat{\mathfrak{sl}(2,\C)}$ and use Theorems \ref{mainintwopthm} and \ref{maintheorem2} to exhibit new intertwining operators among $V^{M(0)}$-modules. We also show by counterexample that the conditions of Theorem \ref{maintheorem} are necessary in order for the conclusions to hold.

First we need to determine explicitly which generalized Verma modules appear in the Garland-Lepowsky resolutions of standard $\widehat{\mathfrak{sl}(2,\C)}$-modules. Recalling notation from Sections 2 and 3, we consider
the action of the subset $\mathcal{W}^1$ of the Weyl group on dominant integral weights of $\widehat{\mathfrak{sl}(2,\mathbb{C})}$, which have the form
\begin{equation*}
 \Lambda=n\dfrac{\alpha}{2}+\ell\mathbf{k'}
\end{equation*}
for $\ell\in\N$ and $0\leq n\leq\ell$. (Note that for the purposes of
studying $\widehat{\mathfrak{sl}(2,\mathbb{C})}$-modules, the coefficient of
$\mathbf{d'}$ in $\Lambda$ is not important.) We recall the action of the
Weyl group generators $r_0$ and $r_1$ on $\mathfrak{H}^*$ from \eqref{sl2Weyl} and note that since $\mathbf{d'}$ is fixed by $\mathcal{W}$, the action of $\mathcal{W}$ on $\mathfrak{H}^*$ projects to an action of $\mathcal{W}$ on
$\mathbb{C}\alpha\oplus\mathbb{C}\mathbf{k'}$. Using \eqref{sl2Weyl}
above, it is easy to prove by induction on $m$:
\begin{lemma}
 Suppose $\Lambda=n\dfrac{\alpha}{2}+\ell\mathbf{k'}$ is a dominant integral
weight of $\widehat{\mathfrak{sl}(2,\mathbb{C})}$. Then for $m\geq 0$, we have
 \begin{equation*}
  (r_0 r_1)^m r_0\left(\Lambda+\rho\right)-\rho =
\left(2((m+1)(\ell+2)-1)-n\right)\dfrac{\alpha}{2}+\ell\mathbf{k'},
 \end{equation*}
and for $m\geq 1$ we have
\begin{equation*}
 (r_0r_1)^m \left(\Lambda+\rho\right)-\rho
=\left(2m(\ell+2)+n\right)\dfrac{\alpha}{2}+\ell\mathbf{k'}.
\end{equation*}
\end{lemma}
Using the fact that $w_j\in \mathcal{W}^1$ is given by $(r_0r_1)^m r_0 $ if $j=2m+1$ is
odd and by $(r_0r_1)^m$ if $j=2m$ is even, we obtain:
\begin{prop}
 If $\Lambda=n\dfrac{\alpha}{2}+\ell\mathbf{k'}$ is a dominant integral weight
of $\widehat{\mathfrak{sl}(2,\mathbb{C})}$, then for $j\geq 0$ we have
 \begin{equation*}
 w_j(\Lambda+\rho)-\rho = \left((\ell+2)j+\dfrac{\ell}{2}(1-(-1)^j)+(-1)^j
n\right)\dfrac{\alpha}{2}+\ell\mathbf{k'}.
 \end{equation*}
\end{prop}
 Thus the resolution of $L(\Lambda)$ given by Theorem \ref{resolution}
is
\begin{equation}\label{explicitres} \cdots \rightarrow
V^{M((\ell+2)j+\ell(1-(-1)^j)/2+(-1)^j n,\,\ell)} \stackrel{d_j}{\rightarrow}
 \cdots
\stackrel{d_2}{\rightarrow} V^{M(2(\ell+1)-n,\,\ell)}
\stackrel{d_1}{\rightarrow} V^{M(\Lambda)}
\stackrel{\Pi_{\Lambda}}{\rightarrow} L(\Lambda) \rightarrow
0,\end{equation}
For $j\in\N$ and $0\leq n\leq\ell$, we use the notation \begin{equation}\label{mjn}
m(j,n)=(\ell+2)j+\dfrac{\ell}{2}(1-(-1)^j)+(-1)^j n.
\end{equation}

We can now apply Theorem \ref{mainintwopthm} to the case that $V^{M(r)}$ appears in
the resolution of a level $\ell$ standard
$\widehat{\mathfrak{sl}(2,\mathbb{C})}$-module given by Theorem
\ref{resolution}. In particular:
\begin{thm}\label{sl2gvmthm}
 Suppose $p,q\geq 0$ and
 $$m(j+1,n),m(j+2,n)\notin\lbrace p+q,p+q-2,\ldots,\vert p-q\vert\rbrace. $$
Then there is a linear isomorphism
$$\mathcal{V}^{V^{M(m(j,n))}}_{V^{M(p)}\,V^{M(q)}}\cong\mathrm{Hom}_{\mathfrak{
sl}(2,\mathbb{C})}(M(p)\otimes M(q),M(m(j,n))).$$
 In particular,
 \begin{equation*}
  \mathcal{N}^{V^{M(m(j,n))}}_{V^{M(p)}\,V^{M(q)}}=\left\lbrace
  \begin{array}{rcl}
  1 & \mathrm{if} & m(j,n)\in\lbrace p+q, p+q-2,\ldots,\vert p-q\vert\rbrace\\
  0 & \mathrm{if} &  m(j,n)\notin\lbrace p+q, p+q-2,\ldots,\vert
p-q\vert\rbrace\\
  \end{array}\right. .
 \end{equation*}
\end{thm}
\pf By Theorem \ref{irred} and the resolution (\ref{explicitres}), $J(m(j,n))$
is irreducible and isomorphic to $L(m(j+1,n))$, and $J(m(j+1,n))$ is irreducible
and isomorphic to $L(m(j+2,n))$. Thus the conclusions follow from Theorem
\ref{mainintwopthm} and the fact that as an
$\mathfrak{sl}(2,\mathbb{C})$-module,
\begin{equation}
 M(p)\otimes M(q)\cong M(p+q)\oplus M(p+q-2)\oplus\ldots\oplus M(\vert p-q\vert)
\end{equation}
for any $p,q\geq 0$. \epfv

\begin{exam}\label{gvmexam}
 Suppose we take $p=1$, $q=m(j,n)$, and $r=m(j,n+1)$ for some $j\geq 0$ and
$0\leq n\leq\ell-1$. Since
 \begin{equation*}
  m(j,n+1)=m(j,n)\pm 1,
 \end{equation*}
there is a one-dimensional space of $\mathfrak{sl}(2,\mathbb{C})$-homomorphisms
from $M(1)\otimes M(m(j,n))$ into $M(m(j,n+1))$. Moreover, one can check
directly from the definition of $m(j,n)$ that
\begin{equation*}
 m(j,n)+1<m(j+1,n+1), m(j+2,n+1).
\end{equation*}
Thus we can conclude from Theorem \ref{sl2gvmthm} that
\begin{equation*}
 \mathcal{N}^{V^{M(m(j,n+1))}}_{V^{M(1)}\,V^{M(m(j,n))}} =1.
\end{equation*}
When $j=n=0$, these are simply the intertwining operators of type
$\binom{V^{M(1)}}{V^{M(1)}\,V^{M(0)}}$, and they are obtained from multiples of
the vertex operator of $V^{M(0)}$ acting on $V^{M(1)}$ by skew-symmetry (see for
example \cite{FHL}). However, when either $j$ or $n$ is positive, these
intertwining operators are new.
\end{exam}

We can also apply Theorem \ref{maintheorem2} to the case that $r=m(j,n)$ and
$q=m(j',n')$ for some $j,j'\geq 0$ and $0\leq n, n'\leq\ell$:
\begin{thm}\label{sl2irredthm}
 Suppose the conditions of Theorem \ref{sl2gvmthm} hold and that in addition
 \begin{equation*}
  m(j'+1,n')\notin\lbrace p+m(j,n), p+m(j,n)-2,\ldots,\vert
p-m(j,n)\vert\rbrace.
 \end{equation*}
Then as vector spaces,
$\mathcal{V}^{V^{M(m(j,n))}}_{V^{M(p)}\,V^{M(m(j',n'))}}\cong\mathcal{V}^{L(m(j,
n))}_{V^{M(p)}\,L(m(j',n'))}$. In particular,
\begin{equation*}
  \mathcal{N}^{L(m(j,n))}_{V^{M(p)}\,L(m(j',n'))}=\left\lbrace
  \begin{array}{rcl}
  1 & \mathrm{if} & m(j,n)\in\lbrace p+m(j',n'), p+m(j',n')-2,\ldots,\vert
p-m(j',n')\vert\rbrace\\
  0 & \mathrm{if} &  m(j,n)\notin\lbrace p+m(j',n'), p+m(j',n')-2,\ldots,\vert
p-m(j',n')\vert\rbrace\\
  \end{array}\right. .
 \end{equation*}
\end{thm}
\pf The theorem follows immediately from Theorems \ref{sl2gvmthm} and
\ref{maintheorem2}. To apply Theorem \ref{maintheorem2}, we recall from Theorem
\ref{irred} and \eqref{explicitres} that $J(m(j',n'))$ is irreducible and thus
generated by its lowest weight space $M(m(j'+1,n'))$. \epfv

\begin{exam}
 As in Example \ref{gvmexam}, we consider $p=1$, $q=m(j,n)$, and $r=m(j,n+1)$
for some $j\geq 0$ and $0\leq n\leq\ell-1$. From Example \ref{gvmexam}, we know
that $\mathcal{N}^{V^{M(m(j,n+1))}}_{V^{M(1)}\,V^{M(m(j,n))}} =1$. It can be
checked directly from the definition \eqref{mjn} of $m(j,n)$ that
 \begin{equation*}
  m(j+1,n)>m(j,n+1)+1.
 \end{equation*}
Thus by Theorem \ref{sl2irredthm},
\begin{equation*}
 \mathcal{N}^{L(m(j,n+1))}_{V^{M(1)}\,L(m(j,n))} =1
\end{equation*}
as well. Note that when $j=n=0$, these intertwining operators can be obtained by
skew-symmetry from intertwining operators of type
$\binom{L(1)}{L(0)\,V^{M(1)}}$. Since $L(1)$ is a module for the simple vertex
operator algebra $L(0)$ when $\ell\geq 1$ (see for example \cite{LL}), a basis
for the one-dimensional space of intertwining operators of type
$\binom{L(1)}{L(0)\,V^{M(1)}}$ is $Y_{L(1)}\circ (1_{L(0)}\otimes\Pi_1)$ where
$\Pi_1$ denotes the canonical projection from $V^{M(1)}$ to $L(1)$. As in
Example \ref{gvmexam}, the intertwining operators discussed in this example are
new when $j$ or $n$ is positive.
\end{exam}

We conclude this section by showing with a counterexample that the conclusions of Theorem \ref{maintheorem} may not hold if the conditions do not hold. Given a level $\ell\in\N$, we take $M(p)=M(q)=M(\ell+1)$ and $M(r)=M(0)=\C\mathbf{1}$. There is a one-dimensional space of $\mathfrak{sl}(2,\C)$-module homomorphisms $M(\ell+1)\otimes M(\ell+1)\rightarrow M(0)$ spanned by a nondegenerate invariant bilinear form on $M(\ell+1)$. From \eqref{explicitres}, the maximal proper submodule of $V^{M(0)}$ is generated by $M(2\ell+2)$. Thus the conditions of Theorem \ref{maintheorem} do not hold because there is a one-dimensional space of $\mathfrak{sl}(2,\C)$-module homomorphisms $M(\ell+1)\otimes M(\ell+1)\rightarrow M(2\ell+2)$.

Now suppose $f: M(\ell+1)\otimes M(\ell+1)\rightarrow M(0)$ is an $\mathfrak{sl}(2,\C)$-module homomorphism and we have maps
\begin{equation*}
 \mathcal{Y}_m: M(\ell+1)\otimes M(\ell+1)\rightarrow V^{M(0)}(-m)
\end{equation*}
for $m\in\Z$ such that $\mathcal{Y}_0=f$ and \eqref{sl2comp1} holds. Since we know (see for example \cite{FZ} or \cite{LL}) that $e(-1)^{\ell+1}\mathbf{1}$ generates the maximal proper submodule $J(0)$, Proposition \ref{radicals} implies
\begin{equation*}
 \langle \mathcal{Y}_{-\ell-1}(u)v, e(-1)^{\ell+1}\mathbf{1}\rangle_M=0
\end{equation*}
for $u,v\in M(\ell+1)$. By \eqref{sl2comp1}, this means
\begin{equation*}
 (-1)^{\ell+1}\langle\mathcal{Y}_0(e(0)^{\ell+1}u)v,\mathbf{1}\rangle_M=(-1)^{\ell+1}\langle f((e^{\ell+1}\cdot u)\otimes v),\mathbf{1}\rangle_M=0.
\end{equation*}
If we take $u=v$ to be a lowest weight vector in $M(\ell+1)$, so that $e^{\ell+1}\cdot u$ is a (non-zero) highest weight vector, we see that $f$ must be zero. Thus we cannot construct maps $\mathcal{Y}_m$ satisfying \eqref{sl2comp1} for non-zero $f$.

Interestingly, we do still have a one-dimensional space of intertwining operators of type $\binom{V^{M(0)}}{V^{M(\ell+1)}\,V^{M(\ell+1)}}$, but the image of these operators is contained in the maximal proper submodule $J(0)$. To construct them, note that Theorem \ref{mainintwopthm} implies that there is a one-dimensional space of intertwining operators of type $\binom{V^{M(2\ell+2)}}{V^{M(\ell+1)}\,V^{M(\ell+1)}}$ induced by the one-dimensional space of $\mathfrak{sl}(2,\C)$-module homomorphisms $M(\ell+1)\otimes M(\ell+1)\rightarrow M(2\ell+2)$. Then if $\mathcal{Y}\in\mathcal{V}^{V^{M(2\ell+2)}}_{V^{M(\ell+1)}\,V^{M(\ell+1)}}$ is non-zero, so is $d_1\circ\mathcal{Y}\in\mathcal{V}^{V^{M(0)}}_{V^{M(\ell+1)}\,V^{M(\ell+1)}}$, where $d_1$ is as in \eqref{explicitres}.

\setcounter{equation}{0}
\section{Appendix: Proof of Theorem \ref{intwopext}}

In this section $\mathfrak{g}$ is an arbitrary finite-dimensional complex simple
Lie algebra, and we fix a level $\ell\in\mathbb{N}$. We will heavily use results from formal calculus, especially delta function properties, in our proof of Theorem \ref{intwopext}. These results will typically be used without comment; see \cite{FLM}, Chapters 2 and 8, and \cite{LL}, Chapter 2, for details on formal calculus.

We now proceed with the proof of Theorem \ref{intwopext}; we recall that
$M(\lambda_1)$ and $M(\lambda_2)$ are irreducible finite-dimensional
$\mathfrak{g}$-modules and $W$ is a $V^{M(0)}$-module. We assume that we have a
map
\begin{equation*}\mathcal{Y}: M(\lambda_1)\otimes M(\lambda_2)\rightarrow
W\lbrace x\rbrace
\end{equation*}
satisfying \eqref{comm1} and \eqref{l0comm}. We need to show that $\mathcal{Y}$
has a unique extension to an intertwining operator of type
$\binom{W}{V^{M(\lambda_1)}\,V^{M(\lambda_2)}}$.

 \allowdisplaybreaks First we extend $\mathcal{Y}$ to a map on
$M(\lambda_1)\otimes V^{M(\lambda_2)}$.  Assuming that $\mathcal{Y}(u,x)v$ has
been defined for
$u\in M(\lambda_1)$ and some $v\in V^{M(\lambda_2)}$, we must define
 \begin{equation}\label{recurs1}
  \mathcal{Y}(u,x)h(-n)v=h(-n)\mathcal{Y}(u,x)v-x^{-n}\mathcal{Y}(h(0)u,x)v
 \end{equation}
for any $h\in\mathfrak{g}$, $n>0$, in order for the commutator formula to hold.
This formula determines a unique well-defined map
\begin{equation*}
 \mathcal{Y}: M(\lambda_1)\otimes V^{M(\lambda_2)}\rightarrow W\lbrace x\rbrace
\end{equation*}
since $V^{M(\lambda_2)}\cong U(\widehat{\mathfrak{g}}_-)\otimes M(\lambda_2)$
linearly.

We show that this $\mathcal{Y}$ is lower truncated. Suppose the lowest
conformal weights of $V^{M(\lambda_1)}$, $V^{M(\lambda_2)}$ and $W$ are
$h_1$, $h_2$, and $h_3$ respectively. (We may assume without loss of generality
that all conformal weights of $W$ are congruent mod $\mathbb{Z}$ to some lowest
weight $h_3$.) Then by (\ref{l0comm}), if $v\in
M(\lambda_2)$,
\begin{equation*}
 \mathcal{Y}(u,x)v=\sum_{m\in\mathbb{Z}}\mathcal{Y}_m(u)v\, x^{-m-h}
\end{equation*}
where $h=h_1+h_2-h_3$, and $\mathcal{Y}_m(u)v\in W(-m)$.  Thus
truncation holds when $v\in M(\lambda_2)$, and if we assume that it holds for
any $u\in M(\lambda_1)$ and some particular $v\in V^{M(\lambda_2)}$, then the
definition (\ref{recurs1}) shows that $\mathcal{Y}(u,x)h(-n)v$ is also lower
truncated. This proves lower truncation since $V^{M(\lambda_2)}$ is generated by
$M(\lambda_2)$ as $\widehat{\mathfrak{g}}_-$-module.

We now verify the commutator formula
\begin{equation}\label{comm2}
 [g(m),\mathcal{Y}(u,x)]v=x^m\mathcal{Y}(g(0)u,x)v
\end{equation}
for any $g\in\mathfrak{g}$ and $m\in\mathbb{Z}$. By the definition this holds
for $m<0$, and for $m\geq 0$ it holds for $v\in M(\lambda_2)$ by
(\ref{comm1}).
Thus it is enough to suppose that (\ref{comm2}) holds for some $v\in
V^{M(\lambda_2)}$ and then show that it also holds for $h(-n)v$ where
$h\in\mathfrak{g}$ and $n>0$:
\begin{eqnarray*}
  && [g(m),\mathcal{Y}(u,x)]h(-n)v  =
g(m)\mathcal{Y}(u,x)h(-n)v-\mathcal{Y}(u,x)g(m)h(-n)v\nonumber\\
 && \hspace{3em} =
g(m)h(-n)\mathcal{Y}(u,x)v-g(m)x^{-n}\mathcal{Y}(h(0)u,x)v-\mathcal{Y}(u,
x)g(m)h(-n)v\nonumber\\
 &&\hspace{3em} = h(-n)g(m)\mathcal{Y}(u,x)v+([g,h](m-n)+m\langle
g,h\rangle\delta_{m-n}\ell)\mathcal{Y}(u,x)v\nonumber\\
 &&\hspace{4em}
-x^{-n}\mathcal{Y}(h(0)u,x)g(m)v-x^{m-n}\mathcal{Y}(g(0)h(0)u,x)v-\mathcal{Y}(u,
x)g(m)h(-n)v\nonumber\\
 &&\hspace{3em} =  h(-n)\mathcal{Y}(u,x)g(m)v+h(-n)x^m\mathcal{Y}(g(0)u,x)v
+\mathcal{Y}(u,x)[g(m),h(-n)]v\nonumber\\
&&\hspace{4em}+x^{m-n}\mathcal{Y}([g(0),h(0)]u,x)v
-x^{-n}\mathcal{Y}(h(0)u,x)g(m)v-x^{m-n}\mathcal{Y}(g(0)h(0)u,x)v\nonumber\\
 && \hspace{4em}-\mathcal{Y}(u,x)g(m)h(-n)v\nonumber\\
 &&\hspace{3em} =
(h(-n)\mathcal{Y}(u,x)-\mathcal{Y}(u,x)h(-n)-x^{-n}\mathcal{Y}(h(0)u,
x))g(m)v\nonumber\\
 &&\hspace{4em}
+h(-n)x^m\mathcal{Y}(g(0)u,x)v-x^{m-n}\mathcal{Y}(h(0)g(0)u,x)v\nonumber\\
 &&\hspace{3em} = x^m\mathcal{Y}(g(0)u,x)h(-n)v.
\end{eqnarray*}
We also need to prove the $L(0)$-commutator formula
\begin{equation*}
 [L(0),\mathcal{Y}(u,x)]v=x\dfrac{d}{dx}\mathcal{Y}(u,x)v+\mathcal{Y}(L(0)u,x)v
\end{equation*}
for any $u\in M(\lambda_1)$ and $v\in V^{M(\lambda_2)}$. The formula holds for
$v\in M(\lambda_2)$ by (\ref{l0comm}), so it is enough to assume that it holds
for some $v\in V^{M(\lambda_2)}$ and then show that it holds for any
$h(-n)v$:
\begin{align}\label{l0comm2}
 \left(L(0)-x\dfrac{d}{dx}\right) & \mathcal{Y}(u,x)h(-n)v  =\left(L(0)-x\dfrac{d}{dx}\right)(h(-n)\mathcal{Y}(u,x)v-x^{-n}\mathcal{Y}(h(0)u,x)v)\nonumber\\
 & =h(-n)\left(L(0)-x\dfrac{d}{dx}\right)\mathcal{Y}(u,x)v-x^{-n}\left(L(0)-x\dfrac{d}{dx}\right)\mathcal{Y}(h(0)u,x)v\nonumber\\
 & \hspace{1em}+[L(0),h(-n)]\mathcal{Y}(u,x)v-nx^{-n}\mathcal{Y}(h(0)u,x)v\nonumber\\
 & =h(-n)(\mathcal{Y}(L(0)u,x)v+\mathcal{Y}(u,x)L(0)v+n\mathcal{Y}(u,x)v)\nonumber\\
 &\hspace{1em}-x^{-n}(\mathcal{Y}(L(0)h(0)u,x)v+\mathcal{Y}(h(0)u,x)L(0)v+n\mathcal{Y}(h(0)u,x)v)\nonumber\\
 & =\mathcal{Y}(L(0)u,x)h(-n)v+\mathcal{Y}(u,x)h(-n)L(0)v+n\mathcal{Y}(u,x)h(-n)v\nonumber\\
 & =\mathcal{Y}(L(0)u,x)h(-n)v+\mathcal{Y}(u,x)L(0)h(-n)v.
\end{align}

We also want to verify that $\mathcal{Y}$ satisfies the following iterate
formula for $u\in M(\lambda_1)$, $v\in V^{M(\lambda_2)}$, $g\in\mathfrak{g}$,
and $m\geq 0$:
\begin{align}\label{assoc1}
 \mathcal{Y}(g(m)u,x) & = \mathrm{Res}_{x_1} ((x_1-x)^m
g(x_1)\mathcal{Y}(u,x)-(-x+x_1)^m\mathcal{Y}(u,x)g(x_1)\nonumber\\
 & = \mathrm{Res}_{x_1} (x_1-x)^m [g(x_1),\mathcal{Y}(u,x)].
\end{align}
Now, since
\begin{align*}
 [g(x_1), & \mathcal{Y}(u,x)]  = \sum_{n\in\mathbb{Z}}
[g(n),\mathcal{Y}(u,x)]x_1^{-n-1}\nonumber\\
 & = \sum_{n\in\mathbb{Z}} x^n x_1^{-n-1}\mathcal{Y}(g(0)u,x)= x_1^{-1}\delta\left(\frac{x}{x_1}\right)\mathcal{Y}(g(0)u,x)
\end{align*}
and
\begin{equation*}
 (x_1-x)^m x_1^{-1}\delta\left(\frac{x}{x_1}\right)=(x-x)^m x_1^{-1}
\delta\left(\frac{x}{x_1}\right)
\end{equation*}
if $m\geq 0$, we have that
\begin{equation*}
\mathrm{Res}_{x_1} (x_1-x)^m
[g(x_1),\mathcal{Y}(u,x)]=\left\lbrace\begin{array}{ccc} \mathcal{Y}(g(0)u,x) &
\mathrm{if} & m=0\\ 0 &\mathrm{if} & m>0\\\end{array} \right. ,
 \end{equation*}
as desired since $g(m)u=0$ for $m>0$.

Now that we have extended $\mathcal{Y}$ to $M(\lambda_1)\otimes
V^{M(\lambda_2)}$,
we want to extend $\mathcal{Y}$ to $V^{M(\lambda_1)}\otimes V^{M(\lambda_2)}$:
assuming
that $\mathcal{Y}(u,x)$ has been defined and is lower truncated for some
particular $u\in V^{M(\lambda_1)}$, we must define
\begin{equation}\label{recurs2}
\mathcal{Y}(h(-n)u,x)=\mathrm{Res}_{x_1}((x_1-x)^{-n}h(x_1)\mathcal{Y}(u,
x)-(-x+x_1)^{-n}\mathcal{Y}(u,x)h(x_1))
\end{equation}
for any $h\in\mathfrak{g}$ and $n>0$, in order for the iterate formula to
hold. This expression is well defined because both $\mathcal{Y}(u,x)$ and
$h(x_1)$ are lower truncated when acting on $V^{M(\lambda_2)}$. Since
$V^{M(\lambda_1)}\cong U(\widehat{\mathfrak{g}}_-)\otimes M(\lambda_1)$
linearly,
this determines $\mathcal{Y}(u,x)$ for any $u\in V^{M(\lambda_1)}$, provided we
can
show that $\mathcal{Y}(h(-n)u,x)$ is also lower truncated. By the lower
truncation of $\mathcal{Y}(u,x)$, it is clear that the first term in
(\ref{recurs2}) is also lower truncated. As for the second term, for any $v\in
V^{M(\lambda_2)}$,
\begin{equation*}
 \mathrm{Res}_{x_1} (-x+x_1)^{-n}\mathcal{Y}(u,x)h(x_1)v  =  \sum_{i\geq 0}
\binom{-n}{i}(-x)^{-n-i}\mathcal{Y}(u,x)h(i)v.
\end{equation*}
Since $h(i)v=0$ for $i$ sufficiently large, the sum is finite and lower
truncation follows from the lower truncation of $\mathcal{Y}(u,x)$.

We first want to prove the $L(0)$-commutator formula (\ref{l0comm}) for any
$u\in V^{M(\lambda_1)}$. Since it holds for $u\in M(\lambda_1)$ by
(\ref{l0comm2}), it is enough to assume it holds for $u\in V^{M(\lambda_1)}$ and
then show that it holds for $g(-n)u$ as well, where $g\in\mathfrak{g}$ and
$n>0$. In the following calculation, we use the commutation formula
\begin{equation*}
 [L(0),g(x)]=x\dfrac{d}{dx} g(x)+g(x)
\end{equation*}
as well as the fact that the residue of a derivative is zero:
\begin{align}\label{l0comm3}
 &\left(L(0)  -x\dfrac{d}{dx}\right)  \mathcal{Y}(g(-n)u,x)\nonumber\\
 &\hspace{1em}=\mathrm{Res}_{x_1}\left(L(0)-x\dfrac{\partial}{\partial x}\right)\left((x_1-x)^{-n}g(x_1)\mathcal{Y}(u,x)-(-x+x_1)^{-n}\mathcal{Y}(u,x)g(x_1)\right)\nonumber\\
 &\hspace{1em}=\mathrm{Res}_{x_1}\left((x_1-x)^{-n}g(x_1)\left(L(0)-x\dfrac{\partial}{\partial x}\right)\mathcal{Y}(u,x)-(-x+x_1)^{-n}\left(L(0)-x\dfrac{\partial}{\partial x}\right)\mathcal{Y}(u,x)g(x_1)\right)\nonumber\\
 &\hspace{2em}-\mathrm{Res}_{x_1} nx\left((x_1-x)^{-n-1}g(x_1)\mathcal{Y}(u,x)-(-x+x_1)^{-n-1}\mathcal{Y}(u,x)g(x_1)\right)\nonumber\\
 &\hspace{2em}+\mathrm{Res}_{x_1}(x_1-x)^{-n}[L(0),g(x_1)]\mathcal{Y}(u,x)\nonumber\\
 &\hspace{1em}=\mathrm{Res}_{x_1}\left((x_1-x)^{-n}g(x_1)\mathcal{Y}(L(0)u,x)-(-x+x_1)^{-n}\mathcal{Y}(L(0)u,x)g(x_1)\right)\nonumber\\
 &\hspace{2em}+\mathrm{Res}_{x_1}\left((x_1-x)^{-n}g(x_1)\mathcal{Y}(u,x)L(0)-(-x+x_1)^{-n}\mathcal{Y}(u,x)L(0)g(x_1)\right)\nonumber\\
 &\hspace{2em}+\mathrm{Res}_{x_1} n\left((x_1-x)^{-n}g(x_1)\mathcal{Y}(u,x_1)-(-x+x_1)^{-n}\mathcal{Y}(u,x)g(x_1)\right)\nonumber\\
 &\hspace{2em}-\mathrm{Res}_{x_1} nx_1  \left((x_1-x)^{-n-1}g(x_1)\mathcal{Y}(u,x_1)-(-x+x_1)^{-n-1}\mathcal{Y}(u,x)g(x_1)\right)\nonumber\\
 &\hspace{2em}+\mathrm{Res}_{x_1}(x_1-x)^{-n}[L(0),g(x_1)]\mathcal{Y}(u,x)\nonumber\\
 &\hspace{1em} =\mathcal{Y}(g(-n)L(0)u,x)+\mathcal{Y}(g(-n),x)L(0)+n\mathcal{Y}(g(-n)u,x)\nonumber\\
 &\hspace{2em}+\mathrm{Res}_{x_1} \left(x_1\dfrac{\partial}{\partial x_1}\left((x_1-x)^{-n}\right)g(x_1)\mathcal{Y}(u,x)-x_1\dfrac{\partial}{\partial x_1}\left((-x+x_1)^{-n}\right)\mathcal{Y}(u,x)g(x_1)\right)\nonumber\\
 &\hspace{2em}+\mathrm{Res}_{x_1}\left((x_1-x)^{-n}[L(0),g(x_1)]\mathcal{Y}(u,x)-(-x+x_1)^{-n}\mathcal{Y}(u,x)[L(0),g(x_1)]\right)\nonumber\\
  &\hspace{1em} =\mathcal{Y}(L(0)g(-n)u,x)+\mathcal{Y}(g(-n),x)L(0)\nonumber\\
 &\hspace{2em}+\mathrm{Res}_{x_1} \left(\dfrac{\partial}{\partial x_1}\left(x_1(x_1-x)^{-n}\right)g(x_1)\mathcal{Y}(u,x)-\dfrac{\partial}{\partial x_1}\left(x_1(-x+x_1)^{-n}\right)\mathcal{Y}(u,x)g(x_1)\right)\nonumber\\
 &\hspace{2em}-\mathrm{Res}_{x_1} \left((x_1-x)^{-n}g(x_1)\mathcal{Y}(u,x)-(-x+x_1)^{-n}\mathcal{Y}(u,x)g(x_1)\right)\nonumber\\
 &\hspace{2em}+\mathrm{Res}_{x_1}\left(x_1(x_1-x)^{-n}\dfrac{\partial}{\partial x_1} g(x_1)\mathcal{Y}(u,x)-x_1(-x+x_1)^{-n}\mathcal{Y}(u,x)\dfrac{\partial}{\partial x_1} g(x_1)\right)\nonumber\\
 &\hspace{2em}+\mathrm{Res}_{x_1}\left((x_1-x)^{-n}g(x_1)\mathcal{Y}(u,x)-(-x+x_1)^{-n}\mathcal{Y}(u,x)g(x_1)\right)\nonumber\\
 &\hspace{1em} =\mathcal{Y}(L(0)g(-n)u,x)+\mathcal{Y}(g(-n),x)L(0).
\end{align}

Now we prove the Jacobi identity
\begin{align*}
 x_0^{-1}\delta\left(\dfrac{x_1-x_2}{x_0}\right) g(x_1)\mathcal{Y}(u,x_2) & -
x_0^{-1}\delta\left(\dfrac{-x_2+x_1}{x_0}\right)\mathcal{Y}(u,x_2)g(x_1)
\nonumber\\
 & =  x_2^{-1}\delta\left(\dfrac{x_1-x_0}{x_2}\right)\mathcal{Y}(g(x_0)u,x_2).
\end{align*}
for any $g\in\mathfrak{g}$ and $u\in V^{M(\lambda_1)}$. Recalling Remark
\ref{commit}, this is equivalent to
proving the iterate formula
\begin{equation*}
 \mathcal{Y}(g(m)u,x_2)=\mathrm{Res}_{x_1} \left((x_1-x_2)^m
g(x_1)\mathcal{Y}(u,x_2)-(-x_2+x_1)^m \mathcal{Y}(u,x_2)g(x_1)\right)
\end{equation*}
for any $m\in\mathbb{Z}$ and the commutator formula
\begin{equation*}
 [g(x_1),\mathcal{Y}(u,x_2)]=\mathrm{Res}_{x_0}
x_2^{-1}\delta\left(\dfrac{x_1-x_0}{x_2}\right)\mathcal{Y}(g(x_0)u,x_2).
\end{equation*}
By (\ref{comm2}), (\ref{assoc1}), and the definition (\ref{recurs2}), the
commutator and iterate formulas hold for $u\in M(\lambda_1)$, so it
suffices
to show that if they hold for $u\in V^{M(\lambda_1)}$, then they also hold for
$h(-n)u$ where $h\in\mathfrak{h}$ and $n>0$. We prove the iterate formula
first, using the assumptions on $u$ and the generating function form of the affine Lie algebra commutation relations (see \cite{FLM}, Proposition 2.3.1):
\begin{align*}
& \mathcal{Y}  (g(m)  h(-n)u,x_2)=
\mathcal{Y}(h(-n)g(m)u,x_2)+\mathcal{Y}(([g,h](m-n)+m\langle
g,h\rangle\ell\delta_{m-n,0})u,x_2)\nn
 & =   \mathrm{Res}_{y_1} ((y_1-x_2)^{-n}
h(y_1)\mathcal{Y}(g(m)u,x_2)-(-x_2+y_1)^{-n}\mathcal{Y}(g(m)u,
x_2)h(y_1))\nn
 &  \hspace{1em}\,+\mathcal{Y}(([g,h](m-n)+m\langle
g,h\rangle\ell\delta_{m-n,0})u,x_2)\nn
 & =  \mathrm{Res}_{y_1}\mathrm{Res}_{x_1} (y_1-x_2)^{-n}\left((x_1-x_2)^m
h(y_1)g(x_1)\mathcal{Y}(u,x_2)-(-x_2+x_1)^m
h(y_1)\mathcal{Y}(u,x_2)g(x_1)\right)\nn
 & \hspace{1em}\, -\mathrm{Res}_{y_1}\mathrm{Res}_{x_1}
(-x_2+y_1)^{-n}\left((x_1-x_2)^m
g(x_1)\mathcal{Y}(u,x_2)h(y_1)-(-x_2+x_1)^m
\mathcal{Y}(u,x_2)g(x_1)h(y_1)\right)\nn
 & \hspace{1em}\, +\mathcal{Y}(([g,h](m-n)+m\langle
g,h\rangle\ell\delta_{m-n,0})u,x_2)\nn
 & =  \mathrm{Res}_{x_1}\mathrm{Res}_{y_1} (x_1-x_2)^m g(x_1)\left((y_1-x_2)^{-n} h(y_1)\mathcal{Y}(u,x_2)-(-x_2+y_1)^{-n}
\mathcal{Y}(u,x_2)h(y_1)\right)\nn
 & \hspace{1em}\, -\mathrm{Res}_{x_1}\mathrm{Res}_{y_1} (x_1-x_2)^m
(y_1-x_2)^{-n}
[g(x_1),h(y_1)]\mathcal{Y}(u,x_2)\nn
 & \hspace{1em}\, -\mathrm{Res}_{x_1}\mathrm{Res}_{y_1} (-x_2+x_1)^m
\left((y_1-x_2)^{-n} h(y_1)\mathcal{Y}(u,x_2)-(-x_2+y_1)^{-n}\mathcal{Y}(u,x_2)h(y_1)\right)g(x_1)\nn
 & \hspace{1em}\, +\mathrm{Res}_{x_1}\mathrm{Res}_{y_1} (-x_2+x_1)^m
(-x_2+y_1)^{-n}
\mathcal{Y}(u,x_2)[g(x_1),h(y_1)]\nn
 & \hspace{1em}\, +\mathcal{Y}(([g,h](m-n)+m\langle
g,h\rangle\ell\delta_{m-n,0})u,x_2)\nn
 & =  \mathrm{Res}_{x_1} ((x_1-x_2)^m
g(x_1)\mathcal{Y}(h(-n)u,x_2)-(-x_2+x_1)^m\mathcal{Y}(h(-n)u,
x_2)g(x_1))\nn
 & \hspace{1em}\, -\mathrm{Res}_{x_1}\mathrm{Res}_{y_1}
(x_1-x_2)^m(y_1-x_2)^{-n}
y_1^{-1}\delta\left(\dfrac{x_1}{y_1}\right)
[g,h](x_1)\mathcal{Y}(u,x_2)\nn
& \hspace{1em}\, -\mathrm{Res}_{x_1}\mathrm{Res}_{y_1}
(x_1-x_2)^m(y_1-x_2)^{-n}\langle
g,h\rangle\ell\dfrac{\partial}{\partial
y_1}\left(y_1^{-1}\delta\left(\dfrac{x_1}{y_1}\right)\right)\mathcal{Y}(u,
x_2)\nn
 & \hspace{1em}\, +\mathrm{Res}_{x_1}\mathrm{Res}_{y_1}
(-x_2+x_1)^m(-x_2+y_1)^{-n}
y_1^{-1}\delta\left(\dfrac{x_1}{y_1}\right)
\mathcal{Y}(u,x_2)[g,h](x_1)\nn
& \hspace{1em}\, +\mathrm{Res}_{x_1}\mathrm{Res}_{y_1}
(-x_2+x_1)^m(-x_2+y_1)^{-n}\langle
g,h\rangle\ell\dfrac{\partial}{\partial
y_1}\left(y_1^{-1}\delta\left(\dfrac{x_1}{y_1}\right)\right)\mathcal{Y}(u,
x_2)\nn
& \hspace{1em}\, +\mathcal{Y}(([g,h](m-n)+m\langle
g,h\rangle\ell\delta_{m-n,0})u,x_2))\nn
 & =  \mathrm{Res}_{x_1} ((x_1-x_2)^m
g(x_1)\mathcal{Y}(h(-n)u,x_2)-(-x_2+x_1)^m\mathcal{Y}(h(-n)u,
x_2)g(x_1))\nn
& \hspace{1em}\, -\mathrm{Res}_{x_1}
((x_1-x_2)^{m-n}[g,h](x_1)\mathcal{Y}(u,x_2)-(-x_2+x_1)^{m-n}\mathcal{Y}(u,x_2)[
g,h](x_1))\nn
& \hspace{1em}\, +\mathrm{Res}_{x_1} (x_1-x_2)^m \dfrac{\partial}{\partial
x_1}(x_1-x_2)^{-n}\langle g,h\rangle\ell\,\mathcal{Y}(u,x_2)\nn
& \hspace{1em}\, -\mathrm{Res}_{x_1} (-x_2+x_1)^m \dfrac{\partial}{\partial x_1}
(-x_2+x_1)^{-n}\langle g,h\rangle\ell\,\mathcal{Y}(u,x_2)\nn
& \hspace{1em}\, +\mathcal{Y}(([g,h](m-n)+m\langle
g,h\rangle\ell\delta_{m-n,0})u,x_2))\nn
& =  \mathrm{Res}_{x_1} ((x_1-x_2)^m
g(x_1)\mathcal{Y}(h(-n)u,x_2)-(-x_2+x_1)^m\mathcal{Y}(h(-n)u,
x_2)g(x_1))\nn
& \hspace{1em}\, -\mathrm{Res}_{x_1} (n(x_1-x_2)^{m-n-1}\langle
g,h\rangle\ell\,\mathcal{Y}(u,x_2))+m\langle
g,h\rangle\delta_{m-n,0}\ell\,\mathcal{Y}(u,x_2)\nn
 & =  \mathrm{Res}_{x_1} ((x_1-x_2)^m
g(x_1)\mathcal{Y}(h(-n)u,x_2)-(-x_2+x_1)^m\mathcal{Y}(h(-n)u,x_2)g(x_1)),
\end{align*}
which proves the iterate formula for $h(-n)u$.

Now we need to prove the commutator formula for $h(-n)u$. We use the fact that
the Jacobi identity holds for $u$ and so
\begin{equation*}
 (x_1-x_2)^{m}
g(x_1)\mathcal{Y}(u,x_2)-(-x_2+x_1)^m\mathcal{Y}(u,x_2)g(x_1)=\mathrm{Res}_{x_0}
x_0^m x_2^{-1}\delta\left(\dfrac{x_1-x_0}{x_2}\right)\mathcal{Y}(g(x_0)u,x_2)
\end{equation*}
for any $g\in\mathfrak{g}$ and $m\in\mathbb{Z}$:
\begin{align*}
& [ g(  x_1 ),  \mathcal{Y}(h(-n)u,x_2)] =  \mathrm{Res}_{y_1} (y_1-x_2)^{-n}
[g(  x_1), h(y_1)\mathcal{Y}(u,x_2)]\nonumber\\
 & \hspace{1em} -\mathrm{Res}_{y_1}
(-x_2+y_1)^{-n}[g(x_1),\mathcal{Y}(u,x_2)h(y_1)]\nonumber\\
 & = \mathrm{Res}_{y_1}
\left((y_1-x_2)^{-n}[g(x_1),h(y_1)]\mathcal{Y}(u,x_2)+
(y_1-x_2)^{-n}h(y_1)[g(x_1),\mathcal{Y}(u,x_2)]\right)\nonumber\\
&  \hspace{1em}-\mathrm{Res}_{y_1}
\left((-x_2+y_1)^{-n}[g(x_1),\mathcal{Y}(u,x_2)]h(y_1) +
(-x_2+y_1)^{-n}\mathcal{Y}(u,x_2)[g(x_1),h(y_1)]\right)\nonumber\\
& =  \mathrm{Res}_{y_1}
(y_1-x_2)^{-n}\left(y_1^{-1}\delta\left(\dfrac{x_1}{y_1}\right)[g,h]
(x_1)+\langle
g,h\rangle\ell\dfrac{\partial}{\partial
y_1}\left(y_1^{-1}\delta\left(\dfrac{x_1}{y_1}\right)\right)\right)\mathcal{Y}(u
,x_2)\nonumber\\
&  \hspace{1em}-\mathrm{Res}_{y_1}
(-x_2+y_1)^{-n}\mathcal{Y}(u,x_2)\left(y_1^{-1}\delta\left(\dfrac{x_1}{y_1}
\right)[g,h
](x_1)+ \langle
g,h\rangle\ell\dfrac{\partial}{\partial
y_1}\left(y_1^{-1}\delta\left(\dfrac{x_1}{y_1}\right)\right)\right)\nonumber\\
&  \hspace{1em}+\mathrm{Res}_{x_0}\mathrm{Res}_{y_1}
x_2^{-1}\delta\left(\dfrac{x_1-x_0}{x_2}\right)(y_1-x_2)^{-n}h(y_1)\mathcal{Y}
(g(x_0)u,x_2)\nonumber\\
& \hspace{1em}-\mathrm{Res}_{x_0}\mathrm{Res}_{y_1}
x_2^{-1}\delta\left(\dfrac{x_1-x_0}{x_2}\right)(-x_2+y_1)^{-n}\mathcal{Y}
(g(x_0)u,x_2) h(y_1)\nonumber\\
& =
(x_1-x_2)^{-n}[g,h](x_1)\mathcal{Y}(u,x_2)-(-x_2+x_1)^{-n}\mathcal{Y}(u,
x_2)[g,h](x_1)\nonumber\\
& \hspace{1em} -\langle g,h\rangle\ell\dfrac{\partial}{\partial
x_1}\left((x_1-x_2)^{-n}-(-x_2+x_1)^{-n}\right)\mathcal{Y}(u,x_2)\nonumber\\
&\hspace{1em}  +\mathrm{Res}_{x_0}
x_2^{-1}\delta\left(\dfrac{x_1-x_0}{x_2}\right)\mathcal{Y}(h(-n)g(x_0)u,
x_2)\nonumber\\
& =  \mathrm{Res}_{x_0}\left( x_0^{-n}
x_2^{-1}\delta\left(\dfrac{x_1-x_0}{x_2}\right)\mathcal{Y}([g,h](x_0)u,
x_2)\right)+n\langle
g,h\rangle\ell\dfrac{(-1)^n}{n!}\dfrac{\partial^n}{\partial
x_1^n}\left(x_2^{-1}\delta\left(\dfrac{x_1}{x_2}\right)\right)\mathcal{Y}(u,
x_2)\nonumber\\
&\hspace{1em}+\mathrm{Res}_{x_0}
x_2^{-1}\delta\left(\dfrac{x_1-x_0}{x_2}\right)\mathcal{Y}(g(x_0)h(-n)u,
x_2)\nonumber\\
& \hspace{1em} -\mathrm{Res}_{x_0}
x_2^{-1}\delta\left(\dfrac{x_1-x_0}{x_2}\right)\sum_{m\in\mathbb{Z}}x_0^{-m-1}
\mathcal{Y}\left(([g,h](m-n)+m\delta_{m-n,0}\langle g,h\rangle\ell)u,x_2\right)\nonumber\\
& =  \mathrm{Res}_{x_0}
x_2^{-1}\delta\left(\dfrac{x_1-x_0}{x_2}\right)\mathcal{Y}(g(x_0)h(-n)u,x_2).
\end{align*}

We have thus proved that the Jacobi identity holds for $g(-1)\mathbf{1}\in
V^{M(0)}$, where $g\in\mathfrak{g}$, and any $v\in V^{M(\lambda_1)}$, and it is
obvious
that the Jacobi identity holds for $\mathbf{1}$ and any $v\in V^{M(\lambda_1)}$.
To
prove the Jacobi identity for any $u\in V^{M(0)}$ and $v\in V^{M(\lambda_1)}$,
it
suffices to assume that it holds for some specific $u\in V^{M(0)}$ and then show
that it holds for $g(-n)u$ for any $g\in\mathfrak{g}$ and $n>0$. Observing
throughout that each formal expression is well-defined, we obtain
\begin{align}\label{fullJacobi}
  x_0^{-1}\delta\left(\dfrac{x_1-x_2}{x_0}\right)
& Y(g(-n)u,x_1)\mathcal{Y}(v,x_2)  -
x_0^{-1}\delta\left(\dfrac{-x_2+x_1}{x_0}\right)\mathcal{Y}(v,x_2)Y(g(-n)u,
x_1) \nonumber\\
& =\mathrm{Res}_{y_1}  (y_1-x_1)^{-n}
x_0^{-1}\delta\left(\dfrac{x_1-x_2}{x_0}\right)g(y_1)Y(u,x_1)\mathcal{Y}(v,
x_2)\nonumber\\
& \hspace{1em}\,-\mathrm{Res}_{y_1}  (-x_1+y_1)^{-n}
x_0^{-1}\delta\left(\dfrac{x_1-x_2}{x_0}\right)Y(u,x_1)g(y_1)\mathcal{Y}(v,
x_2)\nonumber\\
& \hspace{1em}\,-\mathrm{Res}_{y_1}  (y_1-x_1)^{-n}
x_0^{-1}\delta\left(\dfrac{-x_2+x_1}{x_0}\right)\mathcal{Y}(v,x_2)g(y_1)Y(u,
x_1)\nonumber\\
& \hspace{1em}\,+\mathrm{Res}_{y_1}  (-x_1+y_1)^{-n}
x_0^{-1}\delta\left(\dfrac{-x_2+x_1}{x_0}\right)\mathcal{Y}(v,x_2)Y(u,
x_1)g(y_1).
\end{align}
We analyze the first and third terms on the right side first: they become
\begin{align}\label{Jacobi13}
 \mathrm{Res}_{y_1} & (y_1-x_1)^{-n}
x_0^{-1}\delta\left(\dfrac{-x_2+x_1}{x_0}\right) [g(y_1),\mathcal{Y}(u,x_2)]
Y(u,x_1)\nonumber\\
 &\hspace{1em}\,+\mathrm{Res}_{y_1} (y_1-x_1)^{-n}
x_2^{-1}\delta\left(\dfrac{x_1-x_0}{x_2}\right)
g(y_1)\mathcal{Y}(Y(u,x_0)v,x_2)\nonumber\\
 &=\mathrm{Res}_{y_1}\mathrm{Res}_{y_2} (y_1-x_1)^{-n}
x_0^{-1}\delta\left(\dfrac{-x_2+x_1}{x_0}\right)
x_2^{-1}\delta\left(\dfrac{y_1-y_2}{x_2}\right)\mathcal{Y}(g(y_2)v,x_2)Y(u,
x_1)\nonumber\\
 &\hspace{1em}\,+\mathrm{Res}_{y_1}
(y_1-x_2-x_0)^{-n}x_2^{-1}\delta\left(\dfrac{x_1-x_0}{x_2}\right)
g(y_1)\mathcal{Y}(Y(u,x_0)v,x_2)\nonumber\\
 &=\mathrm{Res}_{y_1}\mathrm{Res}_{y_2} (x_2+y_2-x_1)^{-n}
x_0^{-1}\delta\left(\dfrac{-x_2+x_1}{x_0}\right)
x_2^{-1}\delta\left(\dfrac{y_1-y_2}{x_2}\right)\mathcal{Y}(g(y_2)v,x_2)Y(u,
x_1)\nonumber\\
&\hspace{1em}\,+x_2^{-1}\delta\left(\dfrac{x_1-x_0}{x_2}\right)\sum_{i\geq 0}
\binom{-n}{i}(-x_0)^i\,\mathrm{Res}_{y_1}
(y_1-x_2)^{-n-i}g(y_1)\mathcal{Y}(Y(u,x_0)v,x_2)\nonumber\\
&=\mathrm{Res}_{y_2} (-x_0+y_2)^{-n}
x_0^{-1}\delta\left(\dfrac{-x_2+x_1}{x_0}\right)\mathcal{Y}(g(y_2)v,x_2)Y(u,
x_1)\nonumber\\
&\hspace{1em}\,+x_2^{-1}\delta\left(\dfrac{x_1-x_0}{x_2}\right)\sum_{i\geq 0}
\binom{-n}{i}(-x_0)^i\,\mathrm{Res}_{y_1}
(y_1-x_2)^{-n-i}g(y_1)\mathcal{Y}(Y(u,x_0)v,x_2).
\end{align}
We next analyze the second and fourth terms on the right side of
(\ref{fullJacobi}), which become
\begin{align}\label{Jacobi24}
 -\mathrm{Res}_{y_1} & (-x_0-x_2+y_1)^{-n}
x_0^{-1}\delta\left(\dfrac{x_1-x_2}{x_0}\right)Y(u,x_1)g(y_1)\mathcal{Y}(v,
x_2)\nonumber\\
& \hspace{1em}\,+\mathrm{Res}_{y_1}  (-x_1+y_1)^{-n}
x_0^{-1}\delta\left(\dfrac{-x_2+x_1}{x_0}\right)\mathcal{Y}(v,x_2)Y(u,
x_1)g(y_1)\nonumber\\
& =- x_0^{-1}\delta\left(\dfrac{x_1-x_2}{x_0}\right) Y(u,x_1)\sum_{i\geq
0}\binom{-n}{i}(-x_0)^{-n-i}\,\mathrm{Res}_{y_1} (y_1-x_2)^i
g(y_1)\mathcal{Y}(v,x_2)\nonumber\\
& \hspace{1em}\,+\mathrm{Res}_{y_1}  (-x_1+y_1)^{-n}
x_0^{-1}\delta\left(\dfrac{-x_2+x_1}{x_0}\right)\mathcal{Y}(v,x_2)Y(u,
x_1)g(y_1)\nonumber\\
& =- x_0^{-1}\delta\left(\dfrac{x_1-x_2}{x_0}\right) Y(u,x_1)\sum_{i\geq
0}\binom{-n}{i}(-x_0)^{-n-i}\,\mathrm{Res}_{y_1} (-x_2+y_1)^i
\mathcal{Y}(v,x_2)g(y_1)\nonumber\\
& \hspace{1em}\,-x_0^{-1}\delta\left(\dfrac{x_1-x_2}{x_0}\right)\sum_{i\geq
0}\binom{-n}{i}(-x_0)^{-n-i} Y(u,x_1)\mathcal{Y}(g(i)v,x_2)\nonumber\\
& \hspace{1em}\,+\mathrm{Res}_{y_1}  (-x_1+y_1)^{-n}
x_0^{-1}\delta\left(\dfrac{-x_2+x_1}{x_0}\right)\mathcal{Y}(v,x_2)Y(u,
x_1)g(y_1)\nonumber\\
& =-\mathrm{Res}_{y_1} (-x_0-x_2+y_1)^{-n}
x_0^{-1}\delta\left(\dfrac{x_1-x_2}{x_0}\right)
Y(u,x_1)\mathcal{Y}(v,x_2)g(y_1)\nonumber\\
&  \hspace{1em}\,-
x_0^{-1}\delta\left(\dfrac{x_1-x_2}{x_0}\right)\mathrm{Res}_{y_2}
(-x_0+y_2)^{-n} Y(u,x_1)\mathcal{Y}(g(y_2)v,x_2)\nonumber\\
& \hspace{1em}\,+\mathrm{Res}_{y_1}  (-x_1+y_1)^{-n}
x_0^{-1}\delta\left(\dfrac{-x_2+x_1}{x_0}\right)\mathcal{Y}(v,x_2)Y(u,
x_1)g(y_1)\nonumber\\
& =- x_0^{-1}\delta\left(\dfrac{x_1-x_2}{x_0}\right)\mathrm{Res}_{y_2}
(-x_0+y_2)^{-n} Y(u,x_1)\mathcal{Y}(g(y_2)v,x_2)\nonumber\\
& \hspace{1em}\,-\mathrm{Res}_{y_1} (-x_1+y_1)^{-n}
x_2^{-1}\delta\left(\dfrac{x_1-x_0}{x_2}\right)\mathcal{Y}(Y(u,x_0)v,
x_2)g(y_1)\nonumber\\
& =- x_0^{-1}\delta\left(\dfrac{x_1-x_2}{x_0}\right)\mathrm{Res}_{y_2}
(-x_0+y_2)^{-n} Y(u,x_1)\mathcal{Y}(g(y_2)v,x_2)\nonumber\\
& \hspace{1em}\,-\mathrm{Res}_{y_1} (-x_2+y_1-x_0)^{-n}
x_2^{-1}\delta\left(\dfrac{x_1-x_0}{x_2}\right)\mathcal{Y}(Y(u,x_0)v,
x_2)g(y_1)\nonumber\\
& =- x_0^{-1}\delta\left(\dfrac{x_1-x_2}{x_0}\right)\mathrm{Res}_{y_2}
(-x_0+y_2)^{-n} Y(u,x_1)\mathcal{Y}(g(y_2)v,x_2)\nonumber\\
& \hspace{1em}\,- x_2^{-1}\delta\left(\dfrac{x_1-x_0}{x_2}\right)\sum_{i\geq
0}\binom{-n}{i}(-x_0)^i\,\mathrm{Res}_{y_1}
(-x_2+y_1)^{-n-i}\mathcal{Y}(Y(u,x_0)v,x_2)g(y_1).
\end{align}
If we add (\ref{Jacobi13}) and (\ref{Jacobi24}) back together, we get
\begin{align*}
 -x_2^{-1} &
\delta\left(\dfrac{x_1-x_0}{x_2}\right)\mathrm{Res}_{y_2}(-x_0+y_2)^{-n}\mathcal
{Y}(Y(u,x_0)g(y_2)v,x_2)\nonumber\\
& \hspace{1em}\,+ x_2^{-1}\delta\left(\dfrac{x_1-x_0}{x_2}\right)\sum_{i\geq
0}\binom{-n}{i}(-x_0)^i\,\mathrm{Res}_{y_1}(y_1-x_2)^{-n-i}
g(y_1)\mathcal{Y}(Y(u,x_0)v,x_2)\nonumber\\
&\hspace{1em}\,-x_2^{-1}\delta\left(\dfrac{x_1-x_0}{x_2}\right)\sum_{i\geq
0}\binom{-n}{i}(-x_0)^i\,\mathrm{Res}_{y_1}(-x_2+y_1)^{-n-i}\mathcal{Y}(Y(u,
x_0)v,x_2)g(y_1)\nonumber\\
&=-x_2^{-1}\delta\left(\dfrac{x_1-x_0}{x_2}\right)\mathrm{Res}_{y_2}(-x_0+y_2)^{
-n}\mathcal{Y}(Y(u,x_0)g(y_2)v,x_2)\nonumber\\
& \hspace{1em}\,+ x_2^{-1}\delta\left(\dfrac{x_1-x_0}{x_2}\right)\sum_{i\geq
0}\binom{-n}{i}(-x_0)^i \mathcal{Y}(g(-n-i)Y(u,x_0)v,x_2)\nonumber\\
&=-x_2^{-1}\delta\left(\dfrac{x_1-x_0}{x_2}\right)\mathrm{Res}_{y_2}(-x_0+y_2)^{
-n}\mathcal{Y}(Y(u,x_0)g(y_2)v,x_2)\nonumber\\
&
\hspace{1em}\,+x_2^{-1}\delta\left(\dfrac{x_1-x_0}{x_2}\right)\mathrm{Res}_{y_2}
(y_2-x_0)^{-n}\mathcal{Y}(g(y_2)Y(u,x_0)v,x_2)\nonumber\\
&
=x_2^{-1}\delta\left(\dfrac{x_1-x_0}{x_2}\right)\mathcal{Y}(Y(g(-n)u,x_0)v,x_2).
\end{align*}
This completes the proof of the Jacobi identity.

Finally, to prove that $\mathcal{Y}$ is an intertwining operator, we need to
prove the $L(-1)$-derivative property. From the Jacobi identity, or specifically
from the commutator formula, we have for any $u\in V^{M(\lambda_1)}$ that
\begin{equation*}
 [L(0),\mathcal{Y}(u,x)]=x\mathcal{Y}(L(-1)u,x)+\mathcal{Y}(L(0)u,x).
\end{equation*}
On the other hand, we proved in (\ref{l0comm3}) that
\begin{equation*}
 [L(0),\mathcal{Y}(u,x)]=x\dfrac{d}{dx}\mathcal{Y}(u,x)+\mathcal{Y}(L(0)u,x).
\end{equation*}
It follows that
\begin{equation*}
 \mathcal{Y}(L(-1)u,x)=\dfrac{d}{dx}\mathcal{Y}(u,x),
\end{equation*}
as desired. This completes the proof that $\mathcal{Y}$ is an intertwining
operator extending the original $\mathcal{Y}: M(\lambda_1)\otimes
M(\lambda_2)\rightarrow W\lbrace x\rbrace.$ The uniqueness assertion of the
theorem follows from the observation that the extensions of $\mathcal{Y}$ from
$M(\lambda_1)\otimes M(\lambda_2)$ to $M(\lambda_1)\otimes V^{M(\lambda_2)}$ and
then from $M(\lambda_1)\otimes V^{M(\lambda_2)}$ to $V^{M(\lambda_1)}\otimes
V^{M(\lambda_2)}$ are uniquely determined by the commutator and iterate
formulas, respectively.
\begin{rema}
 If we take $\lambda_1=\lambda_2=0$ and $W=V^{M(0)}$, with
$\mathcal{Y}(\mathbf{1},x)\mathbf{1}=\mathbf{1}$, we see that the argument here
proves the Jacobi identity for $V^{M(0)}$.
\end{rema}

\noindent{\small \sc
Beijing International Center for Mathematical Research, Peking University,
Beijing, China 100084}\\
{\em E--mail address}:
\texttt{robertmacrae@math.pku.edu.cn} \\

\noindent {\small \sc Department of Mathematics, University of Notre Dame,
Notre Dame, IN 46556} \\
{\em E--mail address}:
\texttt{jyang7@nd.edu}

\end{document}